\documentclass{gtpart}
\usepackage{amsmath}
\usepackage{amssymb}
\usepackage{amscd}
\usepackage{latexsym}
\usepackage{graphicx} 
\usepackage{subfigure}

\agtart
\title[Pinwheels and nullhomologous surgery] {Pinwheels and nullhomologous surgery\\ on $4$-manifolds with $\mathbf{b^+=1}$}

\author{Ronald Fintushel} 
\givenname{Ronald}
\surname{Fintushel}
\address{Department of Mathematics\\Michigan State University\\
         East Lansing, MI 48824\\USA}
\email{ronfint@math.msu.edu}
\urladdr{http://www.math.msu.edu/~ronfint}

\author{Ronald J. Stern} 
\givenname{Ronald}
\surname{Stern}
\address{Department of Mathematics\\University of California\\
         Irvine,  CA 92697\\USA}
\email{rstern@uci.edu}
\urladdr{http://math.uci.edu/~rstern}


\arxivreference{arXiv:1004.3049}  %%% please supply if 
\arxivpassword{h69dc}   %%% paper is in the arXiv

\theoremstyle{plain}
\newtheorem{thm}{Theorem}
\newtheorem{lem}{Lemma}
\newtheorem{cor}{Corollary}
\newtheorem{prop}{Proposition}

\theoremstyle{definition}

\newtheorem*{conj}{Conjecture}
\newtheorem*{quest}{Question}

\def \CPb {\overline{\mathbf{CP}}{}^{2}}
\def \CP {\mathbf{CP}{}^{2}} 
\def \CPC {\mathbf{CP}{}^{2}\# \,3\,\overline{\mathbf{CP}}{}^{2}}
\def \cpc{\mathbf{CP}{}^{2}\# \,2\,\overline{\mathbf{CP}}{}^{2}}
\def \R {\mathbf{R}}
\def \Z {\mathbf{Z}}
\def \Sig{\Sigma}

\def \SS {S^2\times S^2}
\def \TS  {T^2\times S^2}
\def \vt {\vartheta}
\def \vp {\varphi}

\def \a {\alpha}
\def \b {\beta}
\def \g {\gamma}
\def \d {\delta}
\def \k {\kappa}
\def \lam {\lambda}
\def\L{\Lambda}

\def \O {\Omega}
\def \s {\sigma}
\def \t {\tau}
\def \z {\zeta}

\def \bd {\partial}
\def \x {\times}
\def \- {\smallsetminus}
\def \C {\subset}

\def \e {\epsilon}

\def \hA {\widehat{A}}

\def \sign{{\text{sign}}}

\def \ssw {\text{SW}}
\def \sw {\mathcal{SW}}

\def \wB {\widetilde{B}}

\def\J{\mathcal{J}}

\def\wh{\widehat}
\def \bbF {\mathbf{F}}
\def \TO {T_0\x T_0}
\def \TOst {T_{0,T}\x T_{0,S}}
\def \Ts {T_{0,S}}
\def \Tt {T_{0,T}}
\def \zx {z_{\xi}}

\begin{document}
\begin{abstract}    % type your abstract below
We present a method for finding embedded nullhomologous tori in standard 4-manifolds which can be utilized to change their smooth structure. As an application, we show how to obtain infinite families of simply connected smooth $4$-manifolds with $b^+=1$ and $b^-=2,\dots,7$, via surgery on nullhomologous tori embedded in the standard manifolds $\CP\#\, k\,\CPb$, $k=2,\dots,7$. 
\end{abstract}
\maketitle

\section{Introduction\label{Intro}}

A primary goal of smooth $4$-manifold theory is to understand the classification up to diffeomorphism of $4$-manifolds in a fixed homeomorphism type. Unfortunately, there is not yet a single smooth $4$-manifold for which this has been accomplished. A more modest, but still unachieved, goal is to understand whether or not the following conjecture is true:

\begin{conj} Let $X$ be a simply connected smooth $4$-manifold with $b^+\ge1$. Then there are infinitely many mutually nondiffeomorphic $4$-manifolds homeomorphic to $X$.
\end{conj}

Given a specific $X$, there are several approaches to finding such an infinite family. It is the authors' contention that the most useful and straightforward approach is to produce `exotic' manifolds by surgery on $X$ itself.
In this paper we promote this point of view by producing infinite families of mutually nondiffeomorphic manifolds homeomorphic to $\CP\#\,k\,\CPb$ for $k=2,\dots,7$ by means of surgeries on nullhomologous tori embedded in $\CP\#\,k\,\CPb$. (See Theorem~\ref{Main} and the theorems that lead to it.)

In a previous paper, \cite{Santeria}, we showed that there are manifolds $R_k$ homeomorphic to $\CP\#\,k\,\CPb$, for $k=5,\dots,8$, such that (like $\CP\#\,k\,\CPb$), the manifolds $R_k$ have trivial Seiberg-Witten invariants, and  we saw that there are nullhomologous tori in $R_k$ such that surgery on them gives rise to infinite families of distinct smooth $4$-manifolds homeomorphic to $\CP\#\,k\,\CPb$. However, we were never able to show that the $R_k$ are actually diffeomorphic to the standard manifolds $\CP\#\,k\,\CPb$.

The main new input of this paper is our technique for finding useful nullhomologous tori in standard 4-manifolds in terms of `pinwheel structures'. Roughly speaking, our goal is to cover the given 4-manifold by 4-balls $C_{i}$ meeting along solid tori $S^{1} \times D^{2}$ in their boundary, then ambiently push out pairs of 2-handles from one 4-ball $C_{i}$ into an adjacent ball $C_{j}$  (i.e. adding 1-handles to $C_{j}$) so as to decompose the manifold into components $C'_{j}$ that meet along $S^{1} \times$(punctured torus). We then aim to identify nullhomologous tori in these $C'_{j}$ upon which there are surgeries altering the smooth structure on $X$. 

More specifically, a pinwheel structure is a generalization of the idea of a $k$-fold symplectic sum which was introduced by M. Symington \cite{S}. The basic idea is this: One has a sequence $\{ (X_i; S_i, T_i)\}$ of $4$-manifolds $X_i$ with embedded surfaces $S_i$, $T_i$, which intersect transversely at a single point and with the genus $g(T_i)=g(S_{i+1})$. Let $C_i$ denote the complement in $X_i$ of a regular neighborhood of the configuration $S_i\cup T_i$. We wish to glue the $C_i$'s together so that the normal circle bundle of $T_i$ is identified with the normal circle bundle of $S_{i+1}$. This can't be done unless the sum of the Euler numbers of the two bundles is $0$. However, we can remove a $4$-ball around the intersection point $S_i\cap T_i$ leaving the normal bundles over $S_i\- D^2$ and $T_i\- D^2$ which can be trivialized. It is then possible to glue each $(T_i\- D^2)\x S^1$ to $(S_{i+1}\- D^2)\x S^1$ to obtain a manifold whose boundary is a torus bundle over the circle.
(Each $\bd B^4\cap C_i \cong S^3\- $ (Hopf link) $\cong T^2\x I$.)  If this boundary is $T^3$, then one can glue in $T^2\x D^2$ to obtain a closed $4$-manifold $X$ with a `pinwheel structure'. We use this terminology because the components $C_i$  fan out around a central torus $T_c = T^2\x \{0\}\C T^2\x D^2$ like a pinwheel. 

Our approach is then to find a pinwheel structure on the standard manifolds $X_{k}=\CP\#\,k\,\CPb$ where the surfaces in question are spheres. Then by ambiently pushing out a pair of 2-handles from each $C_{i}$ into $C_{i+1}$, we aim to obtain a new pinwheel structure on $X_{k}$ where the interface surfaces are now tori and the new components $C'_{i}$ contain nullhomologous tori  upon which surgery changes smooth structures on $X_{k}$. 

Here is the idea for identifying these tori in the $C'_{i}$: Suppose a 4-manifold $X$ contains an embedding of $T^2\x D^2 = S^1\x (S^1\x D^2)$, with the central torus $T=T^2\x\{0\}$ representing a nontrivial homology class in $X$. Define the Bing double $B_T$  of the central torus $T=T^2\x\{0\}$ to be $B_T=S^1\x$ \{Bing double of the core of $S^1\x D^2$\}. (Of course this depends on the splitting $T^2=S^1\x S^1$ and a choice of framing.) Then these Bing tori are a pair of nullhomologous tori in $T^2\x D^2$, and hence in $X$. In earlier work, we showed that there are surgeries on these tori which alter the smooth structure of $X$ provided that the Seiberg-Witten invariant of $X$ is nontrivial. The utility of pinwheel structures on a manifold $X$ arises from the realization that there is a manifold $A$ properly embedded $T^2\x D^2$ that contains the Bing tori $B_T$, that embeds in many of the pinwheel components $C'_{i}$ (that do not contain any homologically essential tori), and such that there are surgeries on these tori altering the smooth structure of $X$. 

The key to this surgery construction is to find surgeries on these Bing tori $T_i$ in copies of $B_T$ which give us a symplectic $4$-manifold $Q$ in which the images of the $T_i$ become Lagrangian tori $\L_i$ and for which $b_1(Q)$ is the number of surgeries. We call these surgeries `standard'. The manifold $Q$ is a `model' in the sense of \cite{RevEng}. This puts us in a situation where we can employ our reverse engineering process \cite{RevEng} to surger the $\L_i$ and obtain a family of distinct smooth manifolds $X_n$ homeomorphic to the given manifold $X$. Since the composition of the standard surgery on $T_i$ with the surgery on $\L_i$ is again a surgery performed on $T_i$, the manifolds $X_n$ can all be constructed by surgeries directly on the Bing tori in $X$.

We wish to emphasize that the families of manifolds constructed in this paper are presumably not new (see, for example,  Akhmedov and Park \cite{AP1,AP2}, Akhmedov, Baykur, and Park \cite{ABP}, Baldridge and Kirk \cite{BK}, Fintushel and Stern \cite{DN}, and Fintushel, Park, and Stern \cite{RevEng}, and the references therein), although this is currently unproved. The fact that they can all be obtained via surgery on fixed standard manifolds is new. 

Surprisingly, it is the construction of exotic manifolds by surgering $\CP\#\,k\,\CPb$ with larger $k$ via pinwheel structures that is the most challenging. As our abstract implies, we have not yet been able to accomplish this for $\CP\#\,8\,\CPb$ (but of course we have for $E(1)=\CP\#\,9\,\CPb$). We explain the cause of the difficulties near the end of Section~\ref{More}.

\section{Reverse engineering}\label{RevE}

In this section we review the notion of reverse engineering families of $4$-manifolds as in \cite{RevEng}. The idea here is for any given smooth simply connected 4-manifold X to find a {\it model} Y for X; i.e. Y has the same Euler characteristic and signature as X, and Y has nontrivial Seiberg-Witten invariants, positive $b_{1}$, and essential tori that can be surgered to reduce $b_{1}$. (When $b^+_X=1$, these Seiberg-Witten invariants are the `small-perturbation invariants' --- those corresponding to the same chamber as the solutions of the unperturbed equations. They are defined unambiguously when $b^-\le9$. See {\it{e.g.}} \cite{DN}.) We begin by introducing the relevant notation. Suppose that $T$ is a torus of self-intersection $0$ with tubular neighborhood $N_T$. Let $\a$ and $\b$ be generators of $\pi_1(T^2)$ and  let $S^1_{\a}$ and $S^1_{\b}$ be loops in $T^3=\bd N_T$ homologous in $N_T$ to $\a$ and $\b$ respectively. Let $\mu_T$ denote a meridional circle to $T$ in $X$. By {\it{$p/q$-surgery}} on $T$ with respect to $\b$ we mean 
\begin{gather*} X_{T,\b}(p/q) =(X\- N_T) \cup_{\vp} (S^1\x S^1 \x D^2),\\
\vp: S^1\x S^1 \x \bd D^2 \to \bd(X\- N_T) 
\end{gather*}
where the gluing map 
satisfies $\vp_*([\bd D^2]) = q [S^1_{\b}] + p [\mu_T]$ in $H_1(\bd(X\- N_T);\Z )$. We denote the `core torus' $S^1\x S^1\x\{0\}\C X_{T,\b}(p/q)$ by $T_{p/q}$. 
This notation depends on the given trivialization $\{ S^1_{\a}, S^1_{\b},\mu_T\}$ for $T^3=\bd N_T$. When the curve $S^{1}_{\b}$ is nullhomologous in $X\- N_T$, then $H_1(X_{T,\b}(1/q);\Z )= H_1(X;\Z )$. In addition, when $T$ itself is nullhomologous,  then $H_1(X_{T,\b}(p/q);\Z )= H_1(X;\Z ) \oplus \Z / p\Z$.

If the homology class of $T$ is primitive in $H_2(X;\Z)$ then the meridian $\mu_T$ is nullhomologous in $X\- N_T$. If also $S^1_\b$ represents a nontrivial class in $H_1(X;\R)$, then for any integer $p$, in $X_{T,\b}(p)$ the meridian to $T_p=T_{p/1}$ is $S^1_\b + p\,\mu_T$ which is homologous to $S^1_\b$ in $X\- N_T$. But $S^1_\b$ is not trivial in $H_1(X\- N_T;\R)$. This means that $T_p$ is a nullhomologous torus in $X_{T,\b}(p)$. The meridian $\mu_T$ to $T$ becomes a loop on $\bd N_{T_p}$ and it is nullhomologous in $X_{T,\b}(p)\- N_{T_p} = X\- N_T$ and has a preferred pushoff $S^1_{\mu_T}$ on $\bd N_{T_p}$. Notice that $0$-surgery on $\mu_T$ in $X_{T,\b}(p)$ gives $(X_{T,\b}(p))_{T_p,\mu_T}(0) = X$.

 If $X$ is a symplectic manifold and $\L$ is  any Lagrangian torus, there is a canonical framing, called the Lagrangian framing, of $N_{\L}$. This framing is uniquely determined by the property that pushoffs of 
 $\L$ in this framing remain Lagrangian. If one performs $1/n$ surgeries ($n\in \Z$) with respect to the pushoff in this framing of any curve $\lam$ on $\L$, then $X_{\L,\lam}(1/n)$  is also a symplectic manifold, and the core torus $\L_{1/n}$ is a Lagrangian torus in the resultant manifold. We refer the reader to  Auroux, Donaldson, and Katzarkov \cite{ADK} for a full discussion of this phenomenon, which is referred to there as {\it{Luttinger surgery}}. 

\begin{thm}\label{SW} Let $X$ be a symplectic $4$-manifold which contains $b=b_1(X)$ disjoint Lagrangian tori $\L_i$ which are primitive in $H_2(X;\Z)$. Suppose that each $\L_i$ contains a simple loop $\lam_i$ such that the collection $\{\lam_i\}$  generates $H_1(X;\R)$. Let $X'$ be the symplectic manifold which is the result of $\e_i=\pm 1$ surgery on each $\L_i$ with respect to $\lam_i$ and the Lagrangian framing for $N_{\L_i}$, and let $\L'=(\L_b)_{\e_b}$, which is a Lagrangian torus in $X'$.  Let $X_n = X'_{\L',\mu_{\L'}}(1/n)$ be the result of $1/n$-surgery on $\L'$ with respect to the loop $S^1_{\mu_{\L'}}$ on $\bd N_{\L'}$.
Then among the manifolds $\{ X_n\}$, infinitely many are pairwise nondiffeomorphic.
\end{thm}

Note that $b_1(X')=0$ and since $S^1_{\mu_T}$ is nullhomologous in $X'\- N_{\L'}$, $b_1(X_n)=0$ for all $n$.  We need also to point out that the preferred pushoff $S^1_{\mu_{L'}}$ on $\bd N_{\L'}$ is not the Lagrangian pushoff, and indeed, we shall see that the manifolds $X_n$, $n\ge2$ admit no symplectic structure. 

\begin{proof} (cf. \cite{RevEng}) Let $X''$ denote the symplectic manifold obtained from $X$ via Lagrangian surgeries on $(\L_i,\lam_i)$ with Lagrangian surgery coefficient $\e_i$ for $i=1,\dots, b-1$.
Then $X'= X''_{\L_b,\lam_b}(\e_b)$ and $X''= X'_{\L',\mu_{\L'}}(0)$. Let $k'$ and $k''$ denote the canonical classes of $X'$ and $X''$; so the Seiberg-Witten invariants are  $\ssw_{X'}(k')=\pm 1$ and $\ssw_{X''}(k'')=\pm 1$.  Since $\L_b$ is Lagrangian, it follows from a theorem of Taubes \cite{T} that there is no other basic class whose restriction to $H_2(X''\- N_{\L_b},\bd)$ agrees with that of $k''$, and the same is true for $X'$ and $k'$. Furthermore, $X'\- N_{\L'}=X''\- N_{\L_b}$, and the restriction of $k''$ to $H_2(X''\- N_{\L_b},\bd)$ agrees with that of $k'$. In $H_2(X_n)$,  there is just one class which restricts to a fixed class in  $H_2(X_n\- N_{\L'_{1/n}},\bd)$ because $\L'_{1/n}$ is nullhomologous. It now follows from Morgan, Mrowka, and Szabo \cite{MMS} that there is a basic class $k_n\in H_2(X_n)$ satisfying $\ssw_{X_n}(k_n)=\pm 1\pm n$. Hence the integer invariants
$\max\{ |\ssw_{X_n}(k)|,\; k \ {\text{basic for $X_n$}} \}$
will distinguish an infinite family of pairwise nondiffeomorphic manifolds among the $X_n$.
\end{proof}

\section{Bing doubling tori in $4$-manifolds}\label{Bing}

\begin{figure}[!ht]
\begin{center}
     \subfigure[Bing double]{\includegraphics[scale=.7]{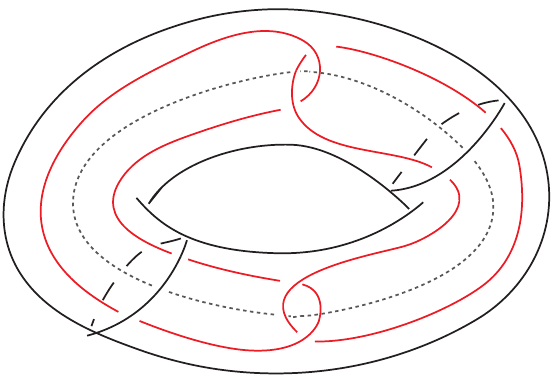}}
     \hspace{.2in}
    \subfigure[Whitehead double]{\includegraphics[scale=.7]{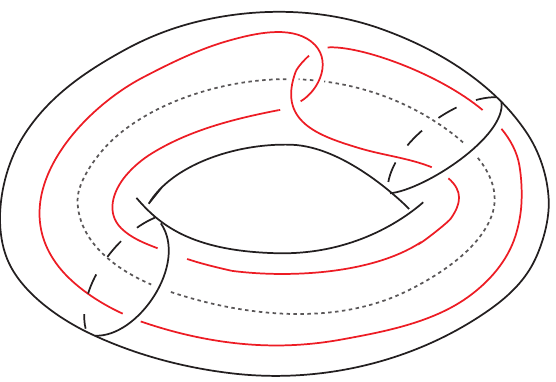}}
      \end{center}
\vspace*{-.1in}     \caption{}
  \label{Doubles}
\end{figure}

Our goal in this section is to see how to find interesting nullhomologous tori in many common $4$-manifolds. As a first step, consider a smooth $4$-manifold $X$ which contains an embedded smooth torus $T$ of self-intersection $0$. Choose local coordinates in which a tubular neighborhood $T\x D^2$ of $T$ is $S^1\x (S^1\x D^2)$. The {\it Bing double} $B_T$ of $T$ consists of the pair of tori $S^1\x$ (Bing double of the core circle $S^1\x \{0\}$). The solid torus $S^1\x D^2$ is shown in Figure~\ref{Doubles}(a). 
This description (including the splitting of $T^2$ into the product $S^1\x S^1$ and a fixed framing, {\it i.e.} a fixed trivialization of the normal bundle of $T$) determines this pair of tori up to isotopy. The component tori in $B_T$ are nullhomologous in $T\x D^2$ and therefore also in $X$.

In order to illustrate the efficacy of surgery on these Bing doubles in changing the smooth structure of $X$, we present a simple example.  Consider the case where $X$ is $E(1)$, the rational elliptic surface, and $T$ is a fiber of a given elliptic fibration on $X$. After splitting off an $S^1$, we choose a framing of the normal bundle of $T$ so that $T\x D^2=S^1\x (S^1\x D^2)$, and each $S^1\x S^1\x \{pt\}$ is an elliptic fiber.
We wish to perform $-1$ surgery on one component of $B_T$ and $1/n$-surgery on the other component to obtain a manifold $X_n$. It is not difficult to see that this manifold is homeomorphic to $X$. The result of the $-1$-surgery on the first component of $B_T$ is to turn the other component into the Whitehead double $\L_W$ of the core torus, namely, $S^1$ times the  Whitehead double of the core circle shown in Figure~\ref{Doubles}(b). Thus $X_n$ is obtained from surgery on $\L_W$ in $X=E(1)$.

We may reinterpret this construction as follows: $X_n$ is obtained by removing a neighborhood $N\cong S^1\x (S^1\x D^2)$  of a fiber torus from $E(1)$, performing $1/n$-surgery on the  Whitehead double torus $\L_W$, and gluing back in. This is the same as forming the fiber sum of $E(1)$ and the result of $1/n$-surgery on $\L_W$ in $S^1\x (S^1\x S^2) = S^1\x (S^1\x D^2\cup S^1\x D^2)$. The particular gluing used for the fiber sum will not matter because the complement of a fiber in $E(1)$ has the property that each diffeomorphism of its boundary extends to a diffeomorphism of its interior. The result of $1/n$-surgery on $\L_W$ in $S^1\x (S^1\x S^2)$ is obtained from $S^1$ times a pair of surgeries on the components of the Whitehead link in $S^3$ with framings $0$ and $1/n$. Performing $1/n$-surgery on one component of the Whitehead link turns the other component into the $n$-twist knot in $S^3$. The upshot of this analysis is that $X_n$ is the fiber sum of $E(1)$ with $S^1$ times $0$-surgery on the $n$-twist knot in $S^3$ --- i.e. $X_n$ is the result of knot surgery on $E(1)$ using the $n$-twist knot. (This shows that $X_n$ is indeed homeomorphic to $X$.) In \cite{KL4M} the (small perturbation) Seiberg-Witten invariant is calculated to be 
$\sw_{X_n}=n\,(t^{-1}-t)$.
Thus surgeries on the Bing double $B_T$ give us an infinite family of distinct smooth $4$-manifolds all homeomorphic to $E(1)$. For a different approach to these same examples, see our paper \cite{Santeria}.

It is the case, however, that there are many $4$-manifolds which do not have any self-intersection $0$ minimal genus tori which one can Bing double; in particular, $\CP\#\,n\,\CPb$, $n\le 8$, and $\SS$. Here, by {\it minimal genus} we mean that there are no spheres that represent the nontrivial homology class of the torus.  We are thus led to ask whether ``Bing doubles" can appear without being embedded in $T^2\x D^2$. Let $T_0= T^2\- {\text{Int}}(D^2)$, a punctured torus. (Here, and throughout this paper, $\text{Int}(S)$ denotes the interior of $S$.) View $B_T\C S^1\x (S^1\x D^2)$, and write the first $S^1$ as $I_1\cup I_2$ ($I_1\cap I_2\cong S^0$), and similarly write the second $S^1$ as $J_1\cup J_2$, and consider $T_0 = S^1\x S^1\- (I_2\x J_2)$. (See Figure~2.) Then $B_T\cap (T_0\x D^2)$ consists of a pair of punctured tori, and $\bd(B_T\cap (T_0\x D^2))$ is a link in $\bd(T_0\x D^2) \cong S^1\!\x S^2\, \# \, S^1\!\x S^2$.

\begin{figure}[ht]
\begin{center}\includegraphics[scale=.8]{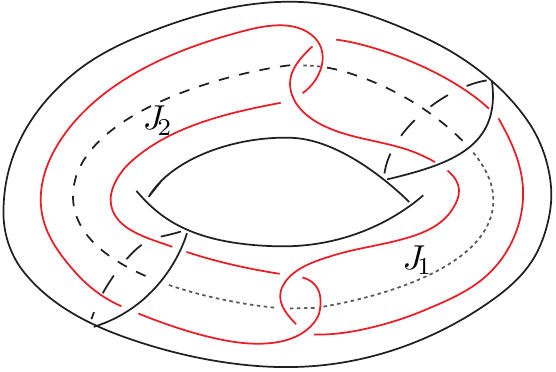}
\end{center}\vspace*{-.1in}
\caption{}
\label{J1J2}
\end{figure}

\begin{figure}[!ht]
\begin{center}
     \subfigure[]{\includegraphics[scale=.3]{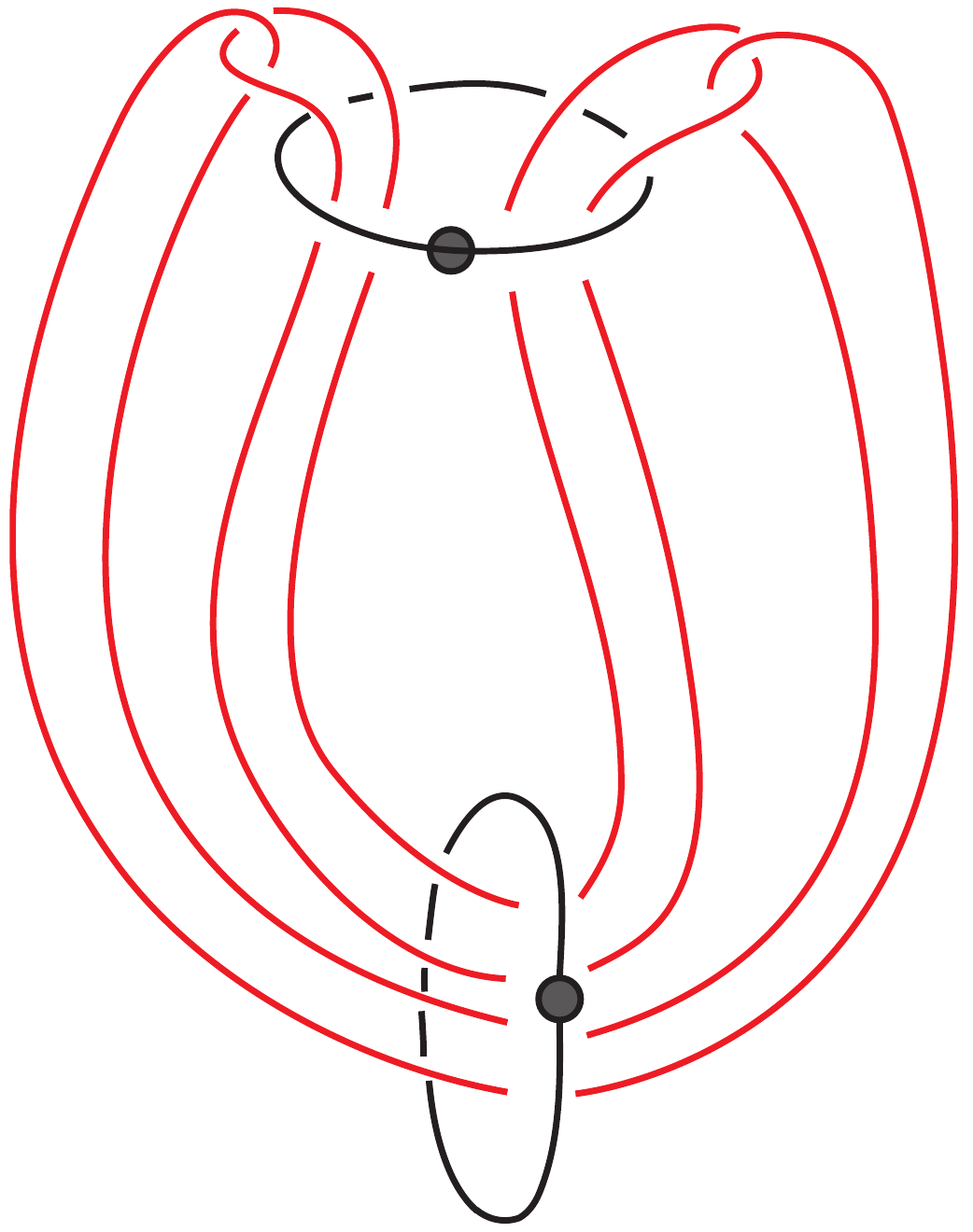}}
     \hspace{.5in}
    \subfigure[]{\includegraphics[scale=.3]{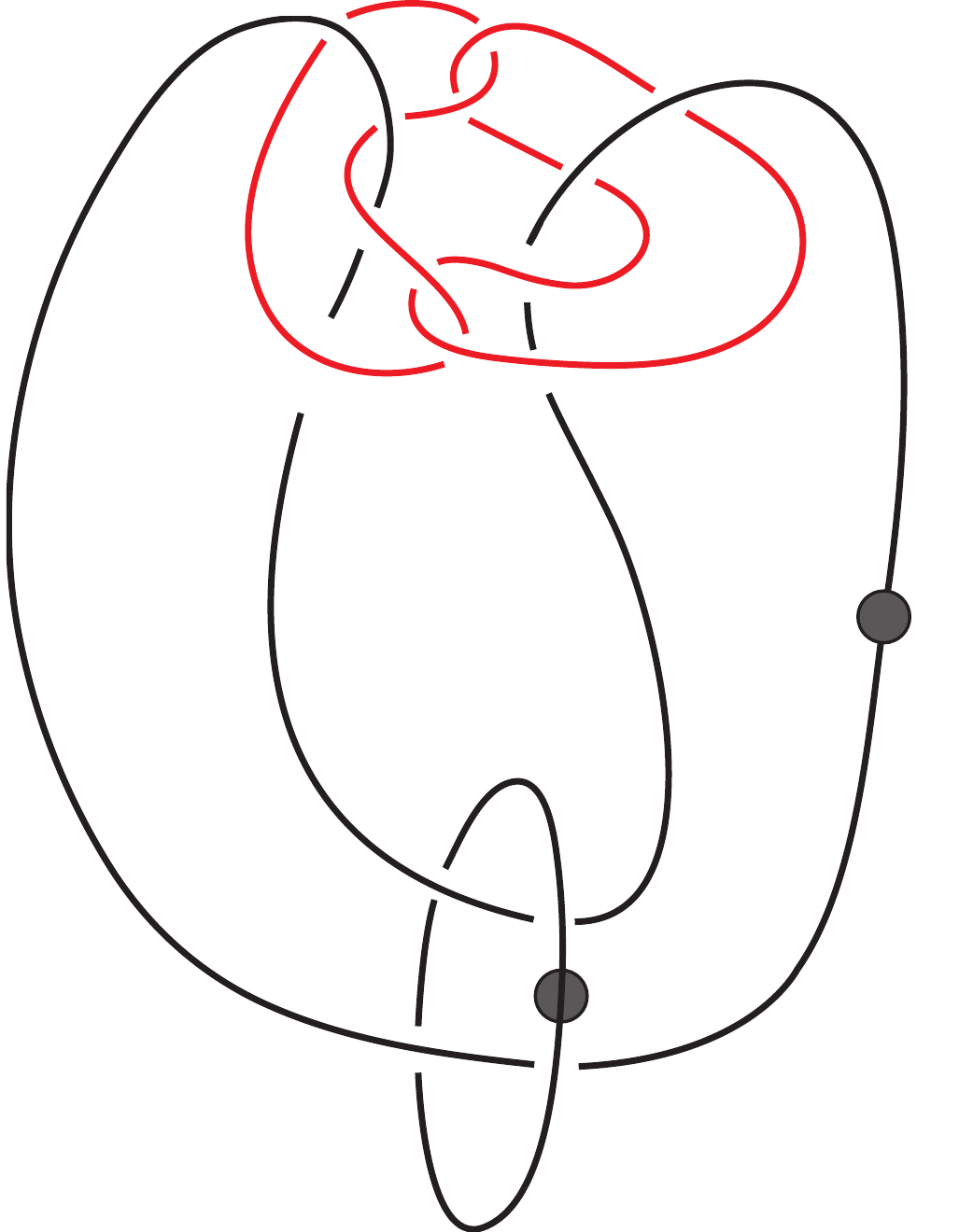}}
      \end{center}
\vspace*{-.2in}     \caption{}
  \label{bdry}
\end{figure}

The intersection of $B_T$ with $(I_2\x J_2)\x D^2$ is a pair of disks $I_2\x$ (intersection of the Bing double link with $J_2\x D^2$). Its boundary is the double of the intersection of the Bing double link with $J_2\x D^2$; i.e. the (1-dimensional) Bing double of $\bd(I_2\x J_2\x \{0\})$. In $\bd(T_0\x D^2) \cong S^1\!\x S^2\, \# \, S^1\!\x S^2$ this boundary is shown in Figure~\ref{bdry}(a) or equivalently, Figure~\ref{bdry}(b).

\begin{figure}[!ht]
\begin{center}
     \subfigure[]{\includegraphics[scale=.3]{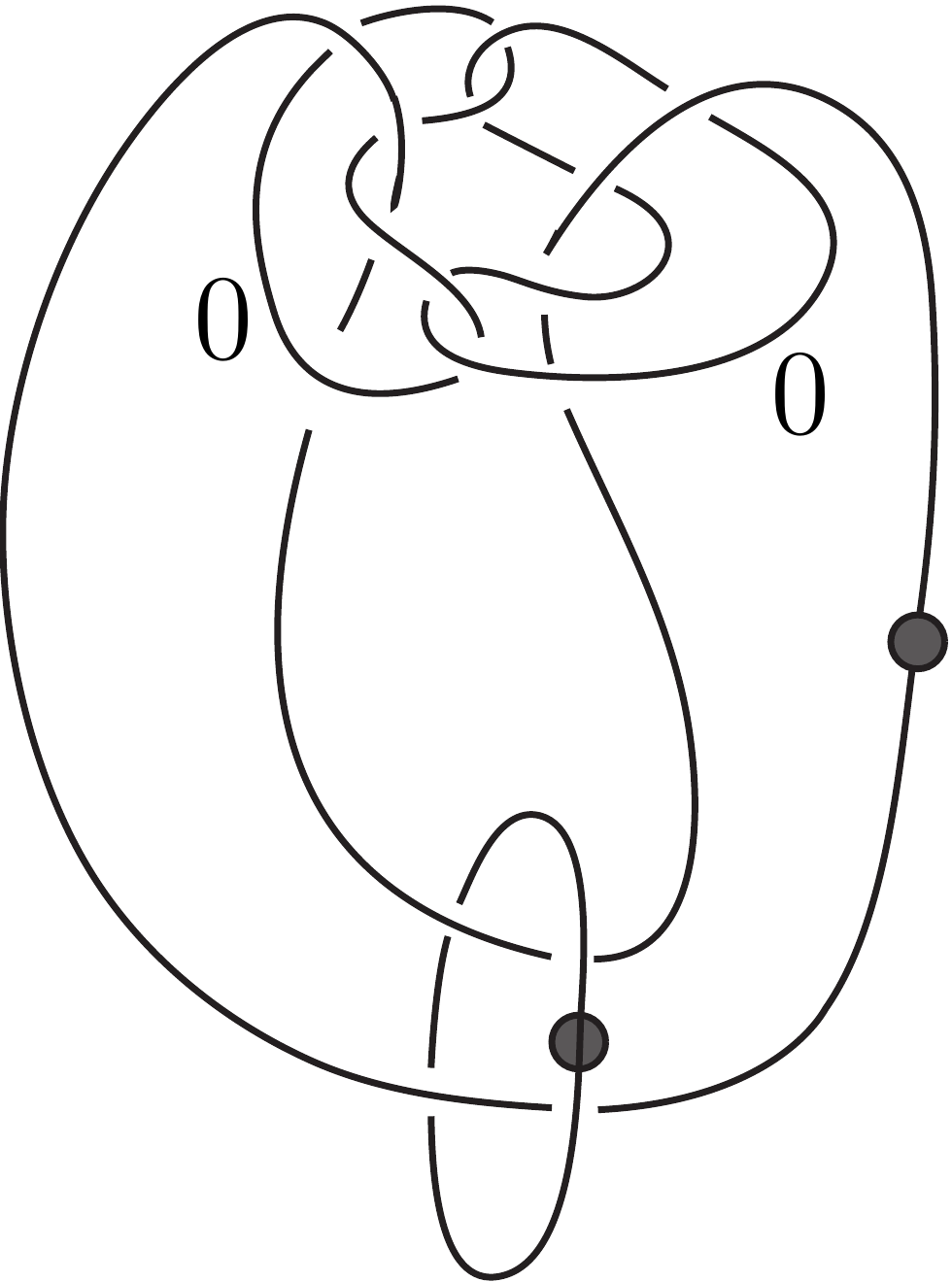}}
     \hspace{.3in}
    \subfigure[]{\includegraphics[scale=1.5]{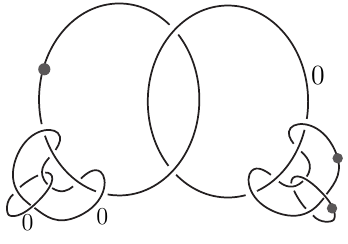}}
      \end{center}
\vspace*{-.1in}     \caption{The manifold $A$}
  \label{A}
\end{figure}

Performing $0$-framed surgeries on these boundary circles (with respect to the framing in Figure~3), we obtain a manifold, shown in Figure~\ref{A}(a) which contains a pair of self-intersection $0$ tori. Call this manifold $A$. It is given equivalently by Figure~\ref{A}(b). (One can see that \ref{A}(a) and \ref{A}(b) are equivalent by twice sliding the bigger Borromean $0$ over the large $0$ in \ref{A}(b) and then cancelling a $1$-$2$ handle pair to obtain \ref{A}(a).) Figure~\ref{A}(b) points out that $A$ is obtained from the $4$-ball by attaching a pair of $2$-handles and then carving out a pair of $2$-handles. The Euler characteristic of $A$ is $e=1$ and its betti numbers are $b_1(A)=2=b_2(A)$.

Our above discussion shows that $A$ embeds in $T^2\x D^2$ and contains the Bing double $B_T$ of the core torus. We continue to refer to the pair of tori in $A$ as $B_T$, even though, in general, they do not constitute a Bing double of some other torus.

\begin{figure}[!ht]
\begin{center}
     \subfigure[]{\includegraphics[scale=.95]{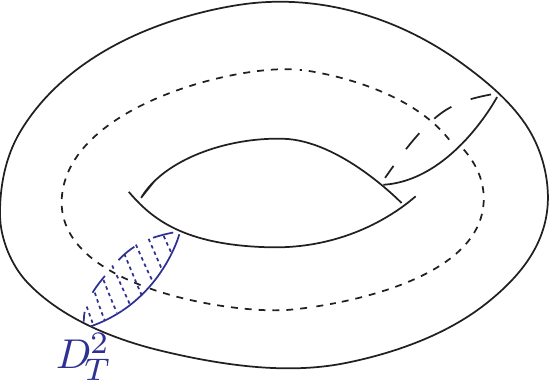}}
     \hspace{.2in}
    \subfigure[]{\includegraphics[scale=.95]{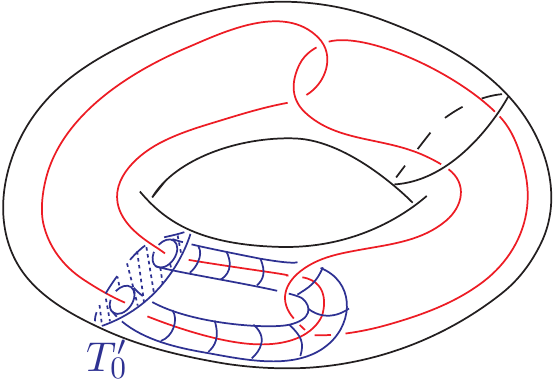}}
      \end{center}
\vspace*{-.1in}     \caption{}
  \label{S^2xT^2}
\end{figure}

\begin{lem}\label{SxT} The manifold $A$ embeds in $T^2\x S^2$ as the complement of a pair of transversely intersecting tori of self-intersection $0$.
\end{lem}
\begin{proof} Write $T^2\x S^2$ as $(S^1\x (S^1\x D^2))\cup (S^1\x (S^1\x D^2))$. We know that $A$ embeds in, say, the second $T^2\x D^2$. In Figure~\ref{A}(a), if we remove the two $2$-handles, we obtain $T_0\x D^2$. The two $2$-handles are attached along the Bing double of the circle $\b=\bd T_0\x\{ 0\}$. If instead, we attach a $0$-framed $2$-handle along $\b$ we obtain $T^2\x D^2$. This implies that $(T^2\x D^2)\- A$ is the complement in the $2$-handle, $D^2\x D^2$, of the core disks of the $2$-handles attached to obtain $A$. This complement is thus the result of attaching two $1$-handles to the $4$-ball. This is precisely the boundary connected sum of two copies of $S^1\x B^3$, {\it i.e.} $T_0\x D^2$. Using the notation in Figure~\ref{J1J2} and above, $(T^2\x D^2)\- (T_0\x D^2)=I_2\x J_2\x D^2$. The complement of the two $2$-handles dug out of this is a neighborhood of $\{pt\}\x$ the shaded punctured torus in Figure~\ref{S^2xT^2}(b). 

Thus the two tori referred to in the lemma are illustrated in Figure~\ref{S^2xT^2}. One of these tori, $T$, is $S^1$ times the core circle in Figure~5(a), and the second torus, $S_T = D_T\cup T_0'$, where $D_T$ is $\{pt\}$ times the shaded meridional disk in Figure~\ref{S^2xT^2}(a), and $T_0'$ is $\{pt\}$ times the shaded punctured torus in Figure~\ref{S^2xT^2}(b). Note that $S_T$ represents the homology class of $\{ pt\}\x S^2$.

From a Kirby calculus point of view, a depiction of a neighborhood $N$ of these two tori is shown in Figure~\ref{NandM}(a) below. Take its union with $A$ as seen in Figure~\ref{A}(b). The Borromean triple on the left side of Figure~\ref{A}(b) cancels with the corresponding triple in Figure~\ref{NandM}(a). We are left with the double of $T^2\x D^2$, i.e. $T^2\x S^2$. 
 \end{proof}

Note that we see $B_T\C S^1\x$ (the solid torus in Figure~\ref{S^2xT^2}(b)) $=T^2\x D^2\C T^2\x S^2$. View $S^1\x S^2$ as $0$-framed surgery on an unknot in $S^3$. The Bing double of the core circle in Figure~\ref{S^2xT^2}(b) is the Bing double of the meridian to the $0$-framed unknot. Performing $0$-framed surgery on the two components of this Bing double gives us $\{0,0,0\}$-surgery on the Borromean rings, {\it viz.} $T^3$. Thus performing $S^1$ times these surgeries gives:

\begin{prop}\label{standard} One can perform surgery on the tori $B_T\C A\C T^2\x S^2$ to obtain the $4$-torus, $T^4$. \qed
\end{prop}

Later we will be interested in other surgeries on $B_T$. We call the surgeries of Proposition~\ref{standard} the {\it standard} surgeries on $B_T$.
Conversely, {\it standard} surgeries on the corresponding pair of tori $\wB_T\C T^4$ yields $T^2\x S^2$. Furthermore, $\wB_T$ is a pair of disjoint Lagrangian tori in $T^4$, $S^1$ times two of the generating circles of $T^3$. The pair of tori of Lemma~\ref{SxT} can also be identified in $T^4$ after the surgeries. The first torus, $S^1$ times the core circle in Figure~5(a) becomes $S^1$ times the third generating circle of $T^3$. Call this torus $T_T$. The other torus intersects $T_T$ once and is disjoint from $B_T$. We call it $T_S$. It is the dual generating torus of $T^4$.  The complement of these tori in $T^4$ is $T_0\x T_0$. We thus have:

\begin{prop}\label{T_0->A} The standard surgeries on the pair of Lagrangian tori $\wB_T$ in $\TO$ give rise to $A$, and conversely, the standard surgeries on $B_T\C A$ yield $\TO$.\qed
\end{prop}

\noindent We wish to emphasize a important point which follows from our discussion.

{\em{The standard surgeries on the Bing tori $B_T$ transform $T^2\x D^2$ into $T^2\x T_0$.}}

Thus the result of the standard surgeries on $T^2\x S^2$ is to transform $A$ into the complement of a transverse pair of generating tori $T_T=T^2\x\{pt\}$ and $T_S=\{pt\}\x T^2$ in $T^4$. The reason for this notation is that $T_S$ is the torus in $T^4$ which is sent to $S_T$ in $T^2\x S^2$ after standard surgeries on $\wB_T$, and $T_T$ is the torus that is sent to $T$.

These surgeries also transform the Bing tori in $B_T$ into the Lagrangian tori $\L_1 = S^1_1\x S^1_3$ and $\L_2 = S^1_1\x S^1_4$ in $T^4=T^2\x T^2=(S^1_1\x S^1_2)\x (S^1_3\x S^1_4)$. The surgeries on $\L_1$ and $\L_2$ are not Lagrangian surgeries in the sense of \cite{ADK}, and so one does not get an induced symplectic structure. Indeed, $T_0\x T_0$ is the complement of transversely intersecting symplectic tori in $T^4$, but after surgery, in $T^2\x S^2$, the complement of $A$ is the regular neighborhood of a pair of tori, one of which is not minimal genus and so cannot be symplectically embedded.

\section{Pinwheels} \label{Pinwheel}

In her thesis (c.f. \cite{S}), Symington discussed the operation of symplectic summing $4$-manifolds along surfaces embedded with normal crossings --- the $k$-fold sum. We study this operation now from a topological point of view. Suppose that we are given smooth $4$-manifolds $X_i$, $i=1,\dots,k$, and that each $X_i$ contains a pair of smoothly embedded surfaces, $S_i$, $T_i$, with genus $g(S_i)$, $g(T_i)$, and self-intersection $m_i$, $n_i$, and suppose that $S_i\cap T_i =\{x_i\}$ is a single transverse intersection. 
Let $N(S_i)$ and $N(T_i)$ be tubular neighborhoods and $N_i$ a $4$-ball neighborhood of $x_i$ large enough so that it is not contained in $N(S_i)\cup N(T_i)$ and such that $S_i$ and $T_i$ intersect it in disks. Then we have 
\[ S_i=S_i'\cup S_i'', \ \ \  T_i=T_i'\cup T_i''\]
where $S_i''$ is a disk, $S_i''=N_i\cap S_i$, and $S_i'=S_i\- S_i''$ is a punctured surface, and similarly for $T_i$. The rough idea of the $k$-fold sum is to remove all the tubular neighborhoods $N(S_i)$ and $N(T_i)$ and glue $\bd N(S_i')$ to $\bd N(T_{i+1}')$, $i=1,\dots,k$ (the subscripts thought of mod $k$), to obtain a manifold whose boundary is a torus bundle over $S^{1}$. Thus we assume that $g(S_i)=g(T_{i+1})$ for all $i$, and our goal is to point out conditions on the self-intersections $m_i$, $n_i$ under which the boundary of the resulting $4$-manifold is the 3-torus $T^{3}$. We then fill in with $T^{2}\times D^{2}$ to obtain a closed $4$-manifold.

To determine the conditions on the self-intersections $m_i$, $n_i$, first note that  the normal bundle of $S_i$ in $X_i$ restricted over the punctured surface  $S_i'$ can be trivialized. Let $\s'_i$ denote the homology class of the normal circle $\{ pt\}\x S^1$ in $H_1(S'_i\x S^1;\Z)$. Similarly let $\t'_i$ denote the homology class of the normal circle in $H_1(T'_i\x S^1;\Z)$. We begin the process of forming the $k$-fold sum by first gluing each $T'_i\x S^1$ to $S'_{i+1}\x S^1$ via a diffeomorphism sending $T'_i$ to $S'_{i+1}$ and $\t_i'$ to $\s_{i+1}'$. Since the boundary of $T'_i$ represents $\s'_i$ and the boundary of $S_{i+1}$ represents $\t'_{i+1}$, in terms of the bases $\{\s'_i, \t'_i\}$ and $\{\s'_{i+1}, \t'_{i+1}\}$, the gluing map restricted to $\bd T'_i\x S^1\to \bd S'_{i+1}\x S^1$ is given by the matrix 
\begin{equation}\label{gluing1} \vp_i=
\left(
\begin{array}{cc}
0 &1 \\
 -1 &0  
\end{array}
\right)
\end{equation}

The tubular neighborhood $N(S_i)= S''_i\x D^2 \cup_{\a_i}  S'_i\x D^2$, where $\a_i: \bd S''_i\x S^1\to  \bd S'_i\x S^1$ is then given by the matrix
\begin{equation}\label{gluing2}  \a_i=
\left(
\begin{array}{cc}
1 & m_i  \\
0 & -1  
\end{array}
\right)
\end{equation}
using bases $\{\s''_i,\t''_i\}$ and $\{\s'_i,\t'_i\}$ where $\s_i''$ and $\t_i''$ are defined in the obvious way.
Also $N(T_i)= T''_i\x D^2 \cup_{\b_i}  T'_i\x D^2$ with
\begin{equation}\label{gluing3}  \b_i=
\left(
\begin{array}{cc}
-1 & 0 \\
n_i & 1  
\end{array}
\right)
\end{equation}
again using bases $\{\s''_i,\t''_i\}$ and $\{\s'_i,\t'_i\}$. 

There is a self-diffeomorphism $r$ of the punctured surface $S_i'$ given by reflection through a plane which restricts to reflection in the boundary circle of $S_i'$. If, instead of gluing $T'_i\x S^1$ to $S'_{i+1}\x S^1$ via the diffeomorphism inducing $\vp_i$ we precompose with $r$ and indentify $\t_i'$ to $-\s_{i+1}'$, we instead get
\[ \vp'_i=
\left(
\begin{array}{cc}
0 &-1 \\
 1 &0  
\end{array}
\right)
\]

\begin{figure}[ht]
\begin{center}
\includegraphics[scale=.9]{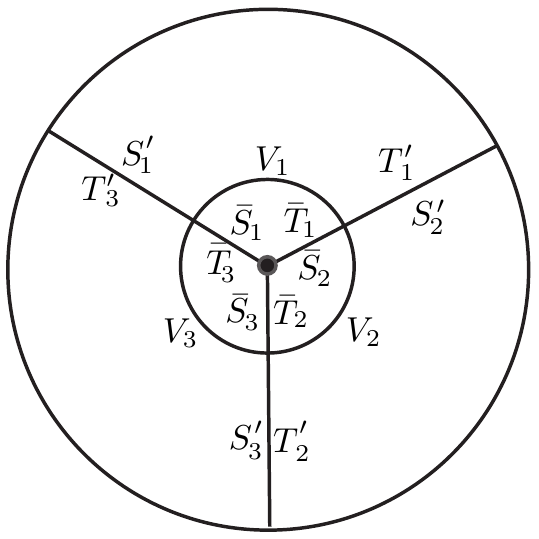}\vspace*{-.1in}
\end{center}
\caption{A 3-fold sum}
\label{3fold}
\end{figure}

Set $V_i=\bd N_i\- (\bd N_i\cap (\text{\em{Int}}(N(S_i)\cup N(T_i)))$. (See Figure~{\ref{3fold}} for a schematic.)  Each $V_i$ is diffeomorphic to the complement of a neighborhood of a Hopf link in $S^3$; i.e. $V_i\cong T^2\x I$. The homology of the factor $T^2$ is generated by $\s_i''$ and $\t_i''$. After gluing together the manifolds $X_i\- (N_i\cup S_i\cup T_i)$, $i=1,\dots,k$, along $T'_i=S'_{i+1}$ as above, we obtain a manifold whose boundary is the union of the $V_i$. This manifold is the torus bundle over the circle whose monodromy is the map $\vt_k\circ\dots\circ\vt_1$ where $\vt_i= \a_{i+1}^{-1}\circ\vp_i\circ\b_i$. Thus
\[ \vt_i=
\left(
\begin{array}{cc}
n_i+m_{i+1}  & 1 \\
-1& 0
\end{array}
\right)
\]
(Compare \cite{S}.) If $\vt_k\circ\dots\circ\vt_1$ is the identity, then the boundary $\bigcup_{i=1}^k V_i\cong T^{2}\times S^{1}$, where $T^{2}$ is the torus fiber generated by the $\s_i''$ and $\t_i''$. This trivial $T^{2}$ bundle over $S^{1}$ can then be extended over the trivial $T^2$ bundle over the disk, $T^2\x D^2$. The manifold thus constructed is called the $k$-fold sum of the $\{ X_i; S_i, T_i\}$. More precisely, the $k$-fold sum $X$ is 
\[ X= \bigcup_{\vp_i}(X_i\- (N\cup N(S_i)\cup N(T_i))) \cup (T^2\x D^2)\]
whenever $\vt_k\circ\dots\circ\vt_1$ is the identity with respect to the given trivializations. If we replace, say, $\vt_k$, with $\vt_k'=
\a_{1}^{-1}\circ\vp'_k\circ\b_k = -\vt_k$, then we see that $X$ can also be constructed when 
$\vt_k\circ\dots\circ\vt_1$ is minus the identity (after changing one of the gluings). In Figure~\ref{3fold}, $\bar{S}_i$ denotes the annulus 
$\bd S'_i \x [0,1]$ which we think of as the disk $S''_i$ with a smaller disk about the origin removed. The $\bar{T}_i$ are defined similarly. We further discuss these annuli below. It is slightly misleading to think of $S_i'$, $\bar{S}_i$, etc. as actually in $X$. Since the normal bundles of these manifolds can be trivialized, they do embed, but in Figure~\ref{3fold}, for example, the text $S_i'$ is only meant to illustrate the boundary of its normal circle bundle, $S_i'\x S^1$, and similarly for other variants of ``S's" and ``T's".

\noindent Note that the tori $T^2\x \{pt\}\C T^2\x D^2$ are nullhomologous in $X$, since $T^2\x \{pt\}$ is isotopic to $\bd(S_i'\x S^1)$ and $\bd(T_i'\x S^1)$.

\begin{prop}{\hspace{-.075in}\rm (Symington \cite{S})}\label{3fd} If $n_i+m_{i+1}=-1$ for $i=1,2,3$ then the $3$-fold sum exists.
\end{prop}
\begin{proof} In this case 
\[ \vt_3\circ\vt_2\circ\vt_1=
\left(
\begin{array}{cc}
-1 & 1 \\
-1& 0 
\end{array}
\right)^3 =
\left(
\begin{array}{cc}
1 & 0 \\
0& 1  
\end{array}
\right)
\]
\end{proof}
In \cite{S}, Symington shows that if all manifolds are symplectic and the surfaces are symplectically embedded then one can arrange for the $3$-fold sum to be symplectic, and similarly in the case of the next proposition.

\begin{prop}{\hspace{-.075in}\rm (McDuff and Symington \cite{MS})}\label{all0}  If $m_{i+1}=-n_i$ for $i=1,\dots,4$ then the $4$-fold sum exists.
\end{prop}
\begin{proof} In this case,
\[ \vt_4\circ\cdots\circ\vt_1=
\left(
\begin{array}{cc}
0 & 1 \\
-1& 0  
\end{array}
\right)^4 =
\left(
\begin{array}{cc}
1 & 0 \\
0& 1  
\end{array}
\right)
\]
\end{proof}

For the next proposition we use the following notation for continued fractions
\[ [c_1,c_2,\dots,c_p] =c_1-1/(c_2-1/(\dots-1/c_p)\dots)\]
In case $c_p=0$, $[c_1,c_2,\dots,c_p]$ is defined to be $[c_1,c_2,\dots,c_{p-2}]$.
For a sequence of integers $\{ a_1,\dots, a_k\}$, we call the set of $k$ continued fractions $[a_1,\dots,a_{k-1}]$, $[a_2,\dots,a_k]$, $[a_3,\dots,a_k,a_1]$, $\dots, [a_k,, a_1\dots,a_{k-2}]$, the {\it cyclic continued fractions} of the sequence. 

\begin{prop} \label{cyclic} Let $a_i=n_i+m_{i+1}$ (indices mod $k$), and consider the sequence $\{ a_1,\dots, a_k\}$. If each of its cyclic continued fractions is $0$, then the corresponding $k$-fold sum exists.
\end{prop}
\begin{proof} We prove this proposition in the case $k=4$, which is the only case that we shall use, and we leave the general case as an exercise for the reader.
\[ \vt_4\circ\cdots\circ\vt_1=
\left(
\begin{array}{cc}
1-a_2a_3 & a_1+a_3-a_1a_2a_3 \\
-a_2-a_4+a_2a_3a_4 & 1-a_1a_2-a_1a_4-a_3a_4 +a_1a_2a_3a_4 
\end{array}
 \right)
\]
The cyclic continued fractions are:
\[ \begin{array}{cc}
(-a_1 - a_3 + a_1a_2 a_3)/(-1 + a_2a_3)& (-a_2 -a_4 + a_2a_3a_4)/(-1 + a_3a_4)\\
 { }& { } \\
(-a_1 - a_3 + a_1 a_3 a_4)/(-1 + a_1 a_4) & (-a_2 - a_4 + a_1 a_2 a_4)/(-1 + a_1 a_2)
\end{array}
\]
The hypothesis immediately implies that the off-diagonal entries of $ \vt_4\circ\cdots\circ\vt_1$ are $0$.
Furthermore, the determinant is $1 - a_1a_2 - a_2a_3 + a_1a_2^2a_3 - a_1a_4 -a_3a_4 +2a_1a_2a_3a_4 + a_2a_3^2a_4 - a_1a_2^2a_3^2a_4$ and it is not difficult to see that the hypothesis shows this $=1$; so $\vt_4\circ\cdots\circ\vt_1 =\pm Id.$
\end{proof}
As mentioned above and shown \cite{S}  and \cite{MS}, if the $(X_i, S_i, T_i)$ are all symplectic, then the conditions in Proposition~ \ref{3fd} and \ref{all0} guarantee that the resulting $k$-fold sum  is symplectic. However, the manifolds obtained as $k$-fold sums using the general criterion of Proposition~\ref{cyclic} need not be symplectic. We give an example below.

\begin{figure}[!ht]
\begin{center}
     \subfigure[Connecting components]{\includegraphics[scale=1.35]{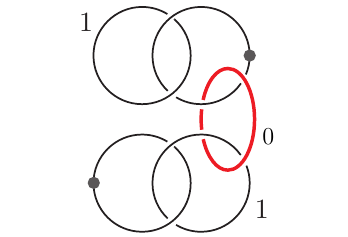}}
     \hspace{.3in} 
    \subfigure[$\CP\smallsetminus T^2\x D^2$]{\includegraphics[scale=1.35]{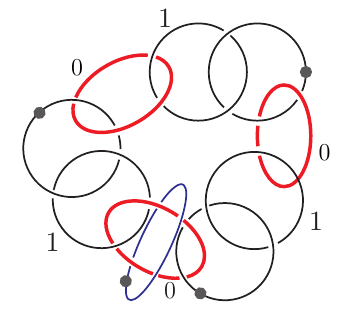}}
      \end{center}
\vspace*{-.1in}     \caption{}
  \label{AroundTheHorn}
\end{figure}

The conditions given in Proposition~\ref{cyclic} are sufficient for the existence of pinwheel structures with given components, but they are not necessary. A necessary and sufficient condition is, of course, that the boundary of the manifold that is formed by gluing the components by identifying the trivial bundles over punctured surfaces in their boundary should be $T^3$. The way to identify these trivial bundles over punctured surfaces in a Kirby calculus diagram is to take the Kirby calculus diagrams of the two components and add a simple loop with framing $0$ around the two meridians that are to be identified. See Figure~\ref{AroundTheHorn}(a).   After all the identifications are made in this way, we introduce a $1$-handle as in  Figure~\ref{AroundTheHorn}(b). This last figure can be seen to have boundary diffeomorphic to $T^3$. It gives an additional check that the pinwheel structure on $\CP$ described in Figure~\ref{CP} exists.

Each $k$-fold sum has a central $2$-torus, $T_c=T^2\x \{0\}\C T^2\x D^2$, and the various component pieces, $X_i\- (\text{\em{Int}}(N(S_i)\cup N(T_i)))$ fan out around $T_c$ like a pinwheel. In order to emphasize this structure, we henceforth refer to a $k$-fold sum as a {\em $k$-component pinwheel} and say that $X$ has a {\em pinwheel structure}. Note that one can have a pinwheel structure without specifying the actual manifolds $X_i$ of a $k$-fold sum. One only needs the complements of neighborhoods of the transversely intersecting surfaces $(S_i, T_i)$ in $X_i$. Of course, given a manifold with the correct boundary, one can form an $X_i$ by gluing in the correct neighborhood of transverse surfaces. We sometimes call $S_i$ and $T_i$ {\em interface surfaces}.

Consider the projection $\pi_i:V_i\cong T^2\x I\to I=I_i$. The central $T^2\x D^2$ extends the trivial $T^2$ bundle over the circle $\bd D^2=\bigcup I_i$. Let $B_i= T^2\x {\text{cone}}(I_i)$ and $C_i=B_i\cup (X_i\- (\text{\em{Int}}(N(S_i)\cup N(T_i))))$. This gives the ``pinwheel structure"
\[ X=\bigcup_{i=1}^k C_i, \ \ \ \ \bigcap _{i=1}^k C_i = T_c \]
We refer to $C_i$ as a {\em{component}} of the pinwheel. Note that $C_i$ is diffeomorphic to 
$X_i\- (\text{\em{Int}}(N(S_i)\cup N(T_i)))$.  We have $\bd S_i'\x S^1\C V_i$, and in fact, if $u_i$ and $v_i$ are the beginning and endpoints of $I_i$, then $\pi_i^{-1}(v_i)= \bd S'_i\x S^1$ 
and $\pi_i^{-1}(u_i)= \bd T'_i\x S^1$. The annuli $\bar{S}_i$ and $\bar{T}_i$ discussed above are the products $\bar{S}_i= [c,v_i] \x \bd S_i\x \{ pt\}$ and  $\bar{T}_i= [c,u_i] \x \bd T_i\x \{ pt\}$, where $[c,v_i]$ and $[c,u_i]$ are the intervals from the cone point to $v_i$ and $u_i$ in ${\text{cone}}(I_i)$. Note that $\pi_i^{-1}[c,v_i]=\bd S_i\x S^1\x [c,v_i]=\bar{S}_i\x S^1$.

Many examples of pinwheels can be obtained from torus actions. Torus actions on simply connected $4$-manifolds are completely classified by Orlik and Raymond \cite{OR} in terms of their orbit data. Briefly, the orbit space must be a disk whose boundary circle is the union of arcs of constant isotropy type separated by (isolated) fixed point images. The isotropy groups which are not trivial or all of $T^2$ are circle subgroups described in polar coordinates by $G(p,q)=\{(\vp,\vt)\mid p\vp+q\vt=0, gcd(p,q)=1\}$.
The orbit space data is a cyclic sequence of pairs $(p_i,q_i)$ of relatively prime integers describing the orbit types over boundary segments. The preimage of each of these closed segments is a $2$-sphere, $A_i$.  Orlik and Raymond's theorem is that each such cyclic sequence gives rise to a simply connected $4$-manifold with a $T^2$-action provided that each determinant 
$\left|\begin{array}{cc}
p_{i-1} & p_i \\
q_{i-1}& q_i  
\end{array}\right| =\pm1$. This condition ensures that the link of each fixed point will be $S^3$. The $2$-sphere $A_i$ which sits over the boundary component of the orbit space which is labelled $(m_i,n_i)$ has self-intersection 
\[ A_i^2=\left|\begin{array}{cc} p_{i-1} & p_i \\q_{i-1}& q_i \end{array}\right| \cdot
\left|\begin{array}{cc} p_i & p_{i+1} \\q_i& q_{i+1} \end{array}\right|\cdot
\left|\begin{array}{cc} p_{i-1}& p_{i+1} \\q_{i-1}& q_{i+1} \end{array}\right| 
=\pm \left|\begin{array}{cc} p_{i-1}& p_{i+1} \\q_{i-1}& q_{i+1} \end{array}\right|  \]
and the intersection number of adjacent spheres is $A_{i-1}\cdot A_i = \left|\begin{array}{cc}
p_{i-1} & p_i \\
q_{i-1}& q_i  
\end{array}\right|$. The second betti number of the $4$-manifold is $b_2= k-2$ where $k$ is the number of fixed points (equivalently the number of distinct segments in the boundary of the orbit space). 
Using this technology, Orlik and Raymond show that a simply connected $4$-manifold which admits a $T^2$-action must be a connected sum of copies of $S^4$, $\CP$, $\CPb$, and $\SS$. So given the orbit data, we may determine the intersection form which tells us which manifold we have.

\noindent{\bf{Example.}} Consider the $T^2$-manifold $X$ with orbit data \[\{(1,-1), (0,1),(1,-1),(2,-1)\}\] The spheres $A_2$ and $A_4$ both have self-intersection $0$ and $A_1^2=-2$, $A_3^2=+2$. Thus $X$ is the ruled surface $\bbF_2 \cong \SS$ with $A_3$ the positive section, $A_1$ the negative section, and $A_2$ and $A_4$ fibers. Here $\bbF_n$ is the ruled surface with positive section having self-intersection $n$.

\begin{figure}[ht]
\vspace*{-.1in}
\begin{center}\includegraphics[scale=1.75]{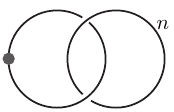}
\end{center}\vspace*{-.1in}
\caption{$B_n$}
\label{Hopf}
\end{figure}

\begin{figure}[!ht]
\begin{center}
     \subfigure[$\CP$]{\includegraphics[scale=.9]{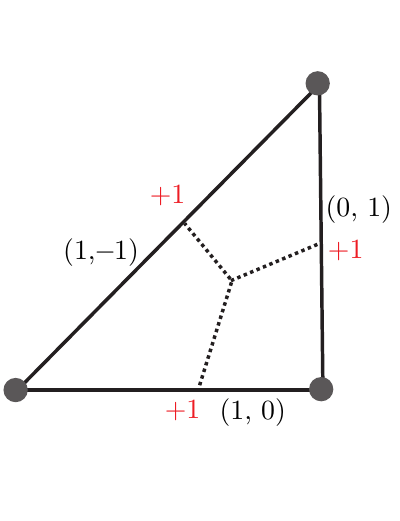}}
    \subfigure[$\CP$ with pinwheel structure]{\includegraphics[scale=.9]{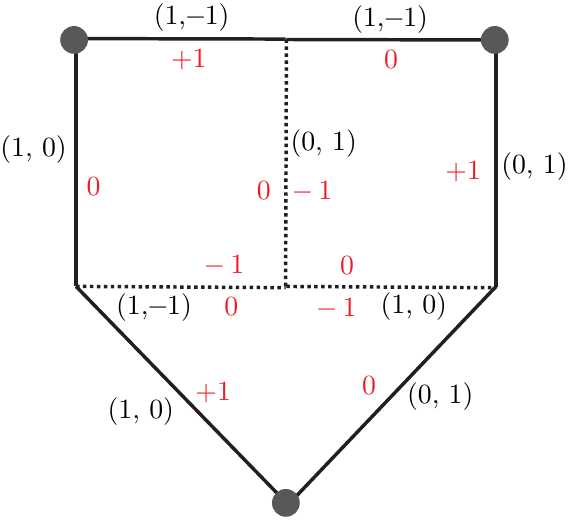}}
     \subfigure[$\CP\#\,\CPb$]{\includegraphics[scale=.9]{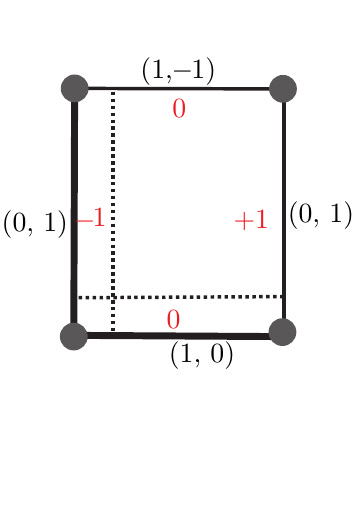}}
      \end{center}
\vspace*{-.1in}     \caption{}
  \label{CP}
\end{figure}

For our next set of examples, note that complement $B_n$ of the negative section and a fiber in the ruled surface $\bbF_n$ is a $4$-ball. This complement in $\bbF_n$ has a handle decomposition with a single cancelling $1$- and $2$-handle pair as shown in Figure~\ref{Hopf}. The framing on the $2$-handle is inherited from the positive section of $\bbF_n$; hence it is $n$ ($\ge 0$).

\smallskip

The next example exhibits $\CP$ as a $3$-component pinwheel made up of three  $4$-balls (the standard coordinate neighborhoods in $\CP$). Our description is essentially that of  Symington \cite{S}, using the language of Orlik and Raymond \cite{OR} rather than that of toric varieties.

\noindent{\bf{Example.}} This example is described by Figure~\ref{CP}. Figure~\ref{CP}(a) shows a torus action on $\CP$, and Figure~\ref{CP}(b) shows how it has the structure of a pinwheel with three $B_{1}$ components.  (See Figure~\ref{CP}(c) where the two removed spheres are indicated by bold lines. This gives the upper right-hand component of Figure~\ref{CP}(b).) The dotted lines in this figure indicate that tubular neighborhoods of the corresponding invariant $2$-spheres have been removed and boundaries glued together. (The point in Figure~\ref{CP}(c) where the two dotted lines meet is the image in the orbit space of the central torus $T_c$.)
The single digits in the figures are self-intersection numbers. Each dotted line represents a solid torus, and the role of the self-intersection numbers assigned to them is to indicate how they get glued together. For example, the boundary of each component in Figure~\ref{CP}(b) (i.e. the union of two solid tori) is $S^3$.
This example illustrates Proposition~\ref{3fd}.

\begin{figure}[!ht]
\begin{center}
     \subfigure[$\CP\#\,\CP$]{\includegraphics[scale=.9]{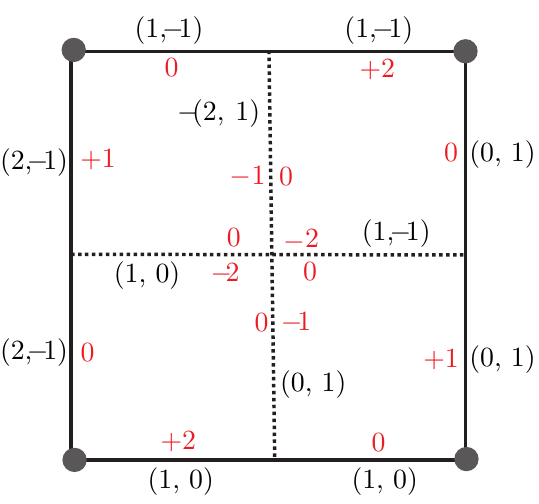}}
     \hspace{.3in} 
    \subfigure[$\SS$]{\includegraphics[scale=.9]{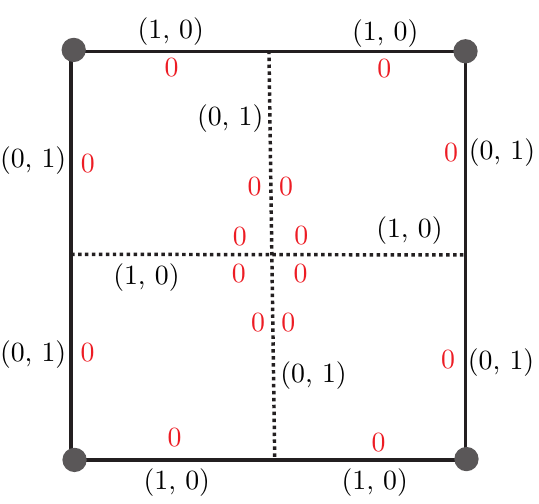}}
      \end{center}
\vspace*{-.1in}     \caption{}
  \label{2CPandSS}
\end{figure}

\noindent{\bf{Example.}} This example is described by Figure~\ref{2CPandSS}(a). It gives $\CP\#\,\CP$ with a $T^2$-action. We use the same notation as in the previous example. This pinwheel structure illustrates Proposition~\ref{cyclic}.  For two of the components we use 
$B_2$ and for the other two, $B_1$.

\noindent{\bf{Example.}} This final example is simpler. It gives $\SS$ as the union of four copies of $B_{0}$, each the complement of  a section and a fiber in $\bbF_0\cong \SS$, and it illustrates Proposition~\ref{all0}. The pinwheel structure is given in Figure~\ref{2CPandSS}(b). 

One can deduce a general technique for constructing pinwheel structures on smooth simply connected manifolds which carry $T^2$-actions: As explained above, the orbit space is a disk which may be viewed as a polygon $P$ whose open edges have a fixed isotropy subgroup  $G(p_i,q_i)$ and whose vertices are the images of fixed points. The pinwheel structure will be built from $n$ ruled surfaces $\bbF_{|r_i|}$ where $r_i$ is the self-intersection number 
\[r_i =  \left|\begin{array}{cc} p_{i-1}& p_{i+1} \\q_{i-1}& q_{i+1} \end{array}\right|\] of the sphere $A_i$ which is the preimage of the $i$th edge under the orbit map. 

\begin{figure}[!ht]
\begin{center}
     \subfigure[]{\includegraphics[scale=1.1]{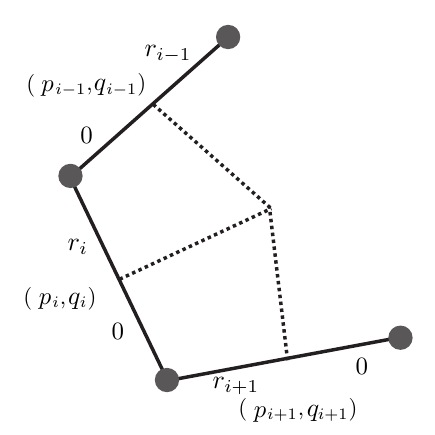}}
     \hspace{.3in} 
    \subfigure[]{\includegraphics[scale=1.1]{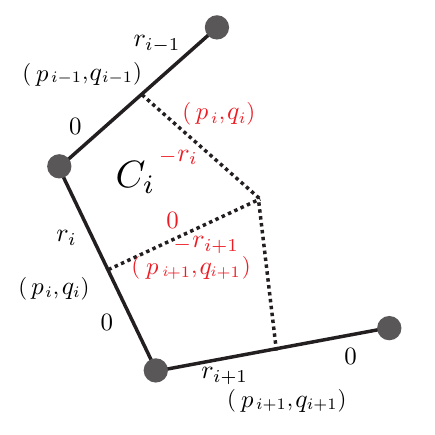}}
      \end{center}
\vspace*{-.1in}     \caption{}
  \label{Barycentric}
\end{figure}

Barycentrically divide $P$ to obtain a union of $4$-gons. Each of the edges will be divided in half and we label the first half-edge with $r_i$ and the second half with $0$. These will be self-intersections of sections or fibers in the ruled surfaces $\bbF_{|r_i|}$. Each (half) edge labelled $0$ must be flanked by two edges with equal isotropy subgroups. (This follows from the rule for computing self-intersections.) Thus in each $4$-gon, two of the isotropy types are given and one is determined by this rule. This determines the isotropy subgroup associated to each edge of the barycentric subdivision. See Figure~\ref{Barycentric} where we are given a $T^2$-action, part of whose polygon is shown in (a), and we deduce the extra data on the barycentric subdivision in (b). The pinwheel component $C_i$ in (b) is obtained from $\bbF_{|r_i|}$ by removing neighborhoods of a fiber and a section of self-intersection $-r_{i}$.

\section{Pinwheel surgery} \label{PS}

As we have seen in \S\ref{Bing}, when we can find an embedding of the manifold $A$ inside a manifold $X$, we get a pair of tori $B_T$ which can be useful to surger. In this section, we show how to accomplish this in some basic $4$-manifolds with pinwheel structures. The general technique we present can be applied in many situations that go beyond the examples given in this paper.

Given a pinwheel structure, the first step is to cyclically push out a pair of $2$-handles from each component into an adjacent component. This has the effect of subtracting a pair of $2$-handles from the adjacent component. The result is that  the genus of the interface increases by one, {\em {e.g.}} a solid torus interface becomes $S^{1} \times T_{0}$, $T_{0}$ a punctured torus.  In other words, we trade handles at the interface of each component. This gives the manifold a new pinwheel structure. The second step is to identify a copy of $A$ inside the new components. The final steps are to compute the effect of the surgeries in $A$ on the Seiberg-Witten invariants and to calculate the fundamental group of the ambient manifold. To compute the effect on Seiberg-Witten invariants, we  shall show that in many cases there are (standard) surgeries on the Bing tori in each of the copies of $A$ that result in a symplectic manifold, but with positive first betti number $b_{1}$. We then employ the techniques of Section~\ref{RevE} to find Lagrangian tori to surger to make $b_{1}=0$. It will follow that there are surgeries on all of the Bing tori in the pinwheel structure that result in distinct smooth manifolds. In order to show that the fundamental group is trivial (so that the resulting manifolds are homeomorphic), we are able to concentrate on the specific pinwheel components  in the examples in Sections~\ref{Pieces} and \ref{Exotic}.

\begin{figure}[ht]
\begin{center}\includegraphics[scale=.75]{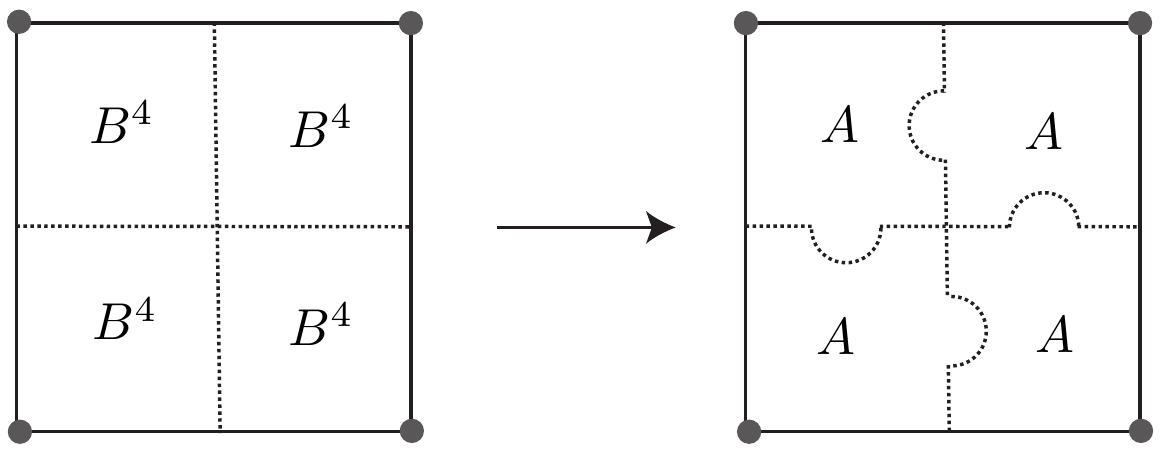}
\end{center}\vspace*{-.1in}
\caption{}
\label{jigsaw}
\end{figure}

The first example to demonstrate this strategy is $\SS$ with the pinwheel structure exhibited in Figure~\ref{2CPandSS}(b). Each component is the complement $B_0$ of a section and a fiber in the ruled surface $\bbF_0$; it is a $4$-ball with  a handle decomposition as given in Figure~\ref{Hopf} with $n=0$. The dotted half-lines in Figure~\ref{2CPandSS}(b) each correspond to a solid torus, and the union of orthogonal half-lines corresponds to the link of a fixed point of the $T^2$-action, a $3$-sphere. This pinwheel structure is represented by the left hand part of Figure~\ref{jigsaw}. We now want to push out $2$-handles from each component into an adjacent component  to obtain a new pinwheel structure represented by the right hand part of Figure~\ref{jigsaw}. To do this, consider the first, say  upper right hand, component  in the left hand part of Figure~\ref{jigsaw}, a $4$-ball  denoted $C_1$. We can view $C_1$ as constructed from a copy of $B^4$ by attaching a $2$-handle with framing $0$ to one component of the Hopf link in Figure~\ref{Hopf} and scooping out a tubular neighborhood of a disk in the $4$-ball bounded by the other component of the Hopf link. (This ``scooping out'' is equivalent to attaching a $1$-handle.) We thus identify the boundary of $C_1$ as the union of two solid tori joined together by a $T^2\x I$; the first solid torus, $U_1$, is the normal circle bundle of the core of the attached $2$-handle, the other solid torus, $V_1$, is the normal circle bundle neighborhood of the scooped out $D^2$, and the $T^2\x I$ is the complement of the Hopf link in $S^3$. We similarly have $\bd C_i= U_i\cup V_i$ for all the components of the pinwheel. We will also use the notation $\bd B_n=U_{(n)}\cup V_{(n)}$, recalling from the previous section that $B_n$ is the complement of the negative section and a fiber in the ruled surface  $\bbF_n$.  The framing on the $2$-handle is inherited from the positive section of $\bbF_n$; hence it is $n$ ($\ge 0$).
The solid torus $U_{(0)}$ has a preferred framing coming from the fact that its core circle is the meridian of the $2$-handle in Figure~\ref{Hopf}, and so it bounds a disk in $\bd B_0$. The Bing double of the core circle of $U_{(0)}$ with respect to this framing bounds disjoint $2$-disks in $B_0$.

The components of the pinwheel structure are glued together so that each $U_i$ is identified with a $V_j$ in an adjoining component. We may assume that the components are ordered so that $V_i$ is identified with $U_{i+1}$.  This induces a preferred framing on $V_i$ with respect to which the Bing double of the core circle of $V_i$ is identified with the Bing double of the core circle of $U_{i+1}$ and so bounds disjoint $2$-disks in $C_{i+1}$.
Add tubular neighborhoods of these disks to $C_i$ and subtract them from $C_{i+1}$. This process has the effect of attaching a pair of $0$-framed $2$-handles to the components of the Bing double link in $V_i\C \bd C_i$ and adding a pair of $1$-handles to the components of the Bing double link in $U_{i+1}\C \bd C_{i+1}$. Once we have applied this process to each pinwheel component $C_i$, we have traded handles so that $\SS=\bigcup_{i=1}^4 \tilde{C_i}$ where $\tilde{C_i}$ is $C_i$ with a pair of $1$-handles and a pair of $2$-handles attached exactly as in Figure~\ref{A}(b). In other words, we have a new pinwheel representation for $\SS$ and each component is a copy of $A$. Figure~\ref{jigsaw} illustrates this process.

\begin{prop} \label{SSsurgery} By performing the standard surgeries on the Bing tori $B_T$ in each of the four copies of $A$ in the above pinwheel structure for $\SS$, one obtains $\Sig_2\x\Sig_2$, the product of two surfaces of genus $2$.\end{prop}
\begin{proof} It follows from Proposition~\ref{T_0->A} that after the standard surgeries on the four pairs of Bing tori, each copy of $A$ is replaced by $T_0\x T_0$. This pinwheel gives $\Sig_2\x\Sig_2$.
\end{proof}

\begin{figure}[!ht]
\begin{center}
     \subfigure[$N$]{\includegraphics[scale=1.3]{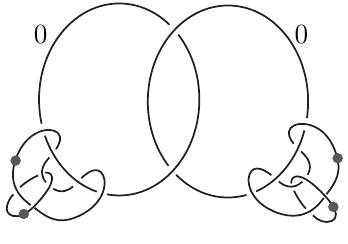}}
     \hspace{.3in} 
    \subfigure[$M$]{\includegraphics[scale=1.3]{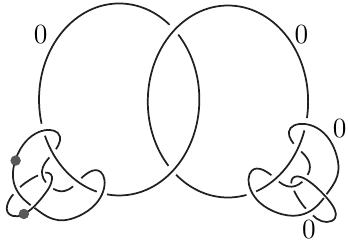}}
      \end{center}
\vspace*{-.2in}     \caption{}
  \label{NandM}
\end{figure}

If we fix one of the copies of $B^4$ in Figure~\ref{jigsaw} then the union of the other three is a regular neighborhood of a section and a fiber in $\bbF_0\cong S^2\x S^2$. After the handle trading process, this neighborhood becomes the manifold $M$ of Figure~\ref{NandM}(b). It is the complement of one copy of $A$ in Figure~\ref{jigsaw}.

Thus $S^2\x S^2=A\cup M$. The manifold $N$ of Figure~\ref{NandM}(a) is the neighborhood of a pair of transversely intersecting tori of self-intersection $0$. Lemma~\ref{SxT} points out that $T^2\x S^2 = A\cup N$, a fact which can also be discerned from the Kirby calculus pictures.

\begin{cor}\label{trade} One can pass from $S^2\x S^2$ to $T^2\x S^2$ by trading a Borromean pair of $2$-handles in $M$ for a Borromean pair of $1$-handles. \qed
\end{cor}
\noindent Note also that the standard surgery on $A$ gives $T_0\x T_0$, and $T^4=(T_0\x T_0) \cup N$.

\medskip

We next formalize the process that we used to obtain the new pinwheel structure for $\SS$. We first need a lemma.

\begin{lem} \label{framing} Consider the framed solid torus $V=D^2\x S^1$ and let $B^{(k)}$ be the $k$-twisted Bing double of the core circle $\{ 0\}\x S^1$. There are unique surgeries on the components of $B^{(k)}$ so that the resulting manifold has $b_1=3$, and the result $W$ of these surgeries is $W=T_0\x S^1$. Let $h\in H_1(\bd(T_0\x S^1))$ be the class of $\{pt\}\x S^1$. Then under the identification  $\bd(T_0\x S^1)= \bd V$,  $h= k[\bd D^2]+[S^1]$. 
\end{lem}
\begin{proof} That there are unique surgeries on $B^{(k)}$ giving a manifold with $b_1=3$ is clear. The union of 
$V\x [0,1]$ with a $2$-handle attached to $\{ 0\}\x S^1\x \{1\}$ with respect to the $k$-twisted framing is a $4$-ball 
$D^2\x D$ where $D$ is the core disk of the $2$-handle. The correspondingly framed Bing double of $\{ 0\}\x S^1\x \{0\}$ bounds a pair of disjoint properly embedded $2$-disks in $D^2\x D$, and their complement is $T_0\x D^2$. The surgeries on $B^{(k)}$ are achieved by this process, {\em{i.e.}} $\bd (T_0\x D^2) = (V\x {0})\cup W$; so $W=T_0\x S^1$.
Because we are using the $k$-framing to attach the $2$-handle to $\{pt\}\x S^1$, it follows that $h= k[\bd D^2]+[S^1]$.
\end{proof}

Now suppose that we have a manifold $X$ with a pinwheel structure with components $C_i$, $i=1,\dots,n$, and suppose also that each $S_i$ (and therefore each $T_i$) is a $2$-sphere. We have 
$ \bd C_i=( T'_i\x S^1)\cup (T^2\x [-1,1]) \cup (S'_i\x S^1)$
as in Section~\ref{Pinwheel}, where $S'_i\cong T'_i\cong D^2$. We suppose that the components are ordered so that $T'_i\x S^1$ is identified with $S'_{i+1}\x S^1$ (via $\vp_i$, see Equation (\ref{gluing1})). 

Consider $\{\text{point}\}\x S^1$ which is a fiber of the normal circle bundle $\bd N(T_i)$. The {\em preferred framing} of this circle is the one it receives from $N(T_i)$. That is, the preferred framing of $\{\text{point}\}\x S^1\C
T_i'\x S^1$ is the one which extends across the normal disk in 
$N(T_i)|_{T_i'}\cong D^2\x D^2$. Of course, the interior of $N(T_i)$ is not contained in $C_i$.
We similarly define the preferred framing for $S'_{i+1}\x S^1$. The condition that we shall need is the:

\smallskip 

\noindent{\bf Handle Trading Condition}.\ Suppose that for each $i$, the circle  $\{pt\}\x S^1\C T'_i\x S^1$ with its preferred framing bounds a proper framed disk $D$ in $C_i$. 

 If this condition holds, then the disk $D$ has a neighborhood $D^2\x D$ in $C_i$ extending the preferred framing on $D^2\x \bd D = T'_i\x S^1$. (For example, consider the case where $C_i$ is obtained from the ruled surface 
$\bbF_n$ by removing tubular neighborhoods of the negative section, $T_i$, and the fiber $S_i$. Then the framed disk $D$ is obtained from that part of another fiber which lies outside of $N(T_i)$.)
As in the proof of Lemma~\ref{framing}, the Bing doubles of $\{pt\}\x S^1$ with respect to this framing bound disjoint disks in $C_i$. We may then use these $2$-disks to attach $2$-handles to $C_{i+1}$ and at the same time a cancelling pair of $1$-handles to $C_i$. We shall refer to this process as {\em handle trading}.
 
 \begin{lem}\label{trading} If a manifold $X$ admits a pinwheel structure with all interface surfaces $S_i$, $T_i$ of genus $0$, and if this pinwheel structure satisfies the the Handle Trading Condition at each interface, then handle trading produces another pinwheel structure on $X$ where all of the interface surfaces are tori.
\end{lem}
\begin{proof} We need to see that the new interface surfaces are tori. Lemma~\ref{framing} implies that the handle trading process turns $T'_i\x S^1$ and $S'_{i+1}\x S^1$ into copies of $T_0\x S^1$ such that the fiber $\{ pt\}\x S^1$ is identified with $\{ pt\}\x S^1$ in both $T'_i\x S^1$ and $S'_{i+1}\x S^1$. (Note that $\vp_i$ preserves the preferred framing.) Hence the new interface surfaces have neighborhoods 
\[ S''_i\x D^2 \cup_{\a_i} T_0\x D^2 \ \ \ {\text {and}} \ \ \  T''_i\x D^2 \cup_{\b_i} T_0\x D^2 \]
so we indeed get a pinwheel structure with tori as interface surfaces.
\end{proof}

\begin{figure}[!ht]
\begin{center}
     \subfigure[ $A_{(1)} $]{\includegraphics[scale=1.2]{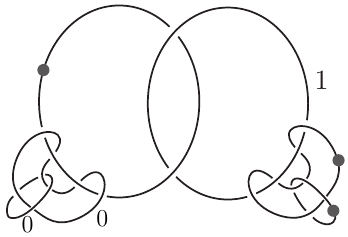}}
     \hspace{.05in} 
    \subfigure[$\hA$]{\includegraphics[scale=1.2]{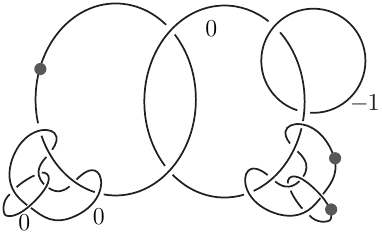}}
    \hspace{.05in} 
    \subfigure[$\hA$]{\includegraphics[scale=.275]{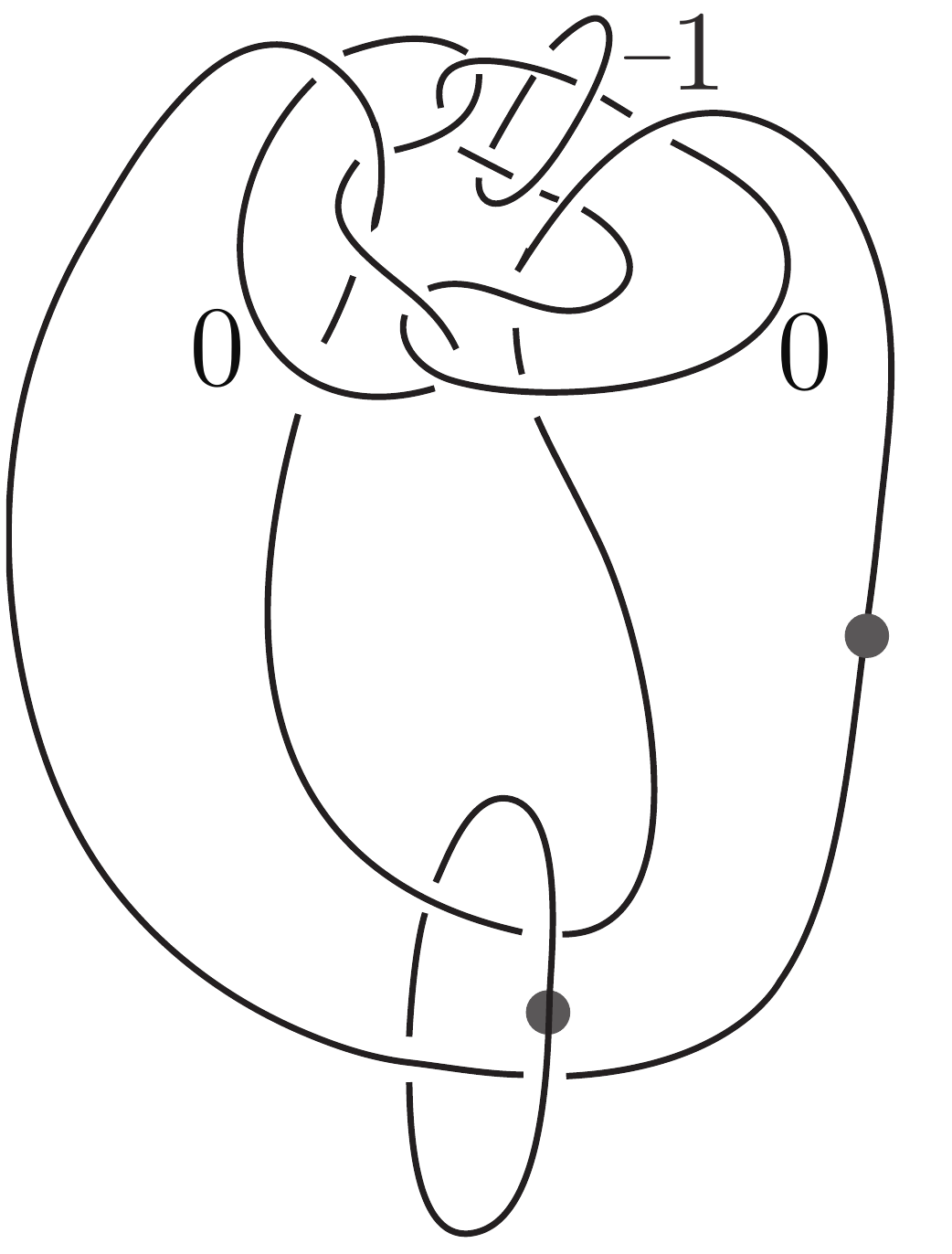}}
      \end{center}
\vspace*{-.1in}     \caption{}
  \label{AA}
\end{figure}

For our next example, consider the pinwheel structure on $\CP$ which is illustrated in Figure~\ref{CP}. Each pinwheel component is given by a pair of cancelling handles as in Figure~\ref{Hopf} with $n=1$. The Handle Trading Condition is satisfied at each interface because the meridian of the `dot' in Figure~\ref{Hopf} is identified with the meridian of a `$+1$' in the next component, and this clearly bounds a disk and the framing extends.
We may again trade handles to obtain the manifold $A_{(1)}$ of Figure~\ref{AA}(a).
If we blow up and slide the $+1$-framed handle in Figure~\ref{AA}(a) over the exceptional curve, we obtain the manifold $\hA$ of Figure~\ref{AA}(b) (which is diffeomorphic to $A_{(1)}\#\,\CPb$), and we see $A\C \hA$. From Figure~\ref{AA}(c) it is clear that at least one of the Bing tori in $A\C \hA$ must intersect the exceptional curve.

 This means that after blowing up in each component of the pinwheel structure of $\CP$, we obtain a $3$-fold pinwheel structure for $\CPC$ on which we can perform our handle trading, and so that  each component of the new pinwheel structure is a copy of $\hA$.  Note that the Euler characteristic of $\hA$ is $e=2$; its betti numbers are $b_1(\hA)=2$, $b_2(\hA)=3$.

For $g\ge 1$, let $\bbF_n(g)$ denote the irrational ruled surface whose base has genus $g$ and with $c_1=n$. Then $A_{(1)} $ is the complement in  $\bbF_1(1)$ of the negative section and another torus representing the homology class of the fiber, but (nonminimal) genus one, similar to Figure~\ref{S^2xT^2}, except now rather than two copies of $T^2\x D^2$ we replace (a) by a tubular neighborhood of the negative section, a $D^2$-bundle over $T^2$ with $c_1=-1$, and (b) with the bundle with $c_1=+1$, a tubular neighborhood of a (positive) section. 

\begin{figure}[!ht]
\begin{center}
     \subfigure[ $M_1$]{\includegraphics[scale=1.35]{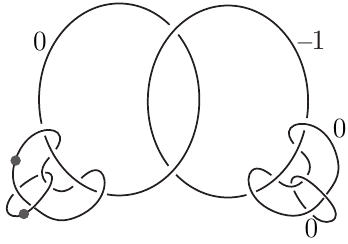}}
     \hspace{.3in} 
    \subfigure[$N_1$]{\includegraphics[scale=1.35]{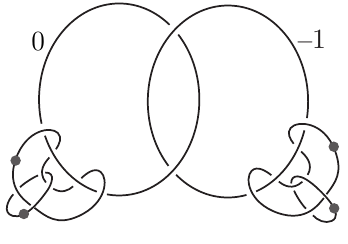}}
      \end{center}
\vspace*{-.1in}     \caption{}
  \label{M1N1}
\end{figure}

The blowup $\bbF_1(1)\#\,\CPb$ is diffeomorphic to 
$\bbF_0(1)\#\,\CPb\cong (T^2\x S^2)\#\,\CPb$. One way to see this is to blow up, separating the fiber $F$ and the negative section $S_-$, and then blow down the new exceptional curve $E'=F-E$ to get 
$\bbF_0(1)$. This correspondence takes $S_-$ to $S'-E'$ where $S'$ is a section of  $\bbF_0(1)$, and it takes fiber to fiber. We have 
\[(\TS)\#\, \CPb = ((T^2\x D^2)\#\, \CPb)\, \cup \,(T^2\x D^2)\]
and to get $\hA$, we remove the proper transform torus $(T^2\x \{0\}) - E'$ and the disk 
$\{ pt\}\x D^2$ from the first summand together with the punctured torus of Figure~\ref{S^2xT^2}(b). We can also see this by an argument similar to that of Corollary~\ref{trade}: Start with 
$\bbF_1\cong \CP\#\,\CPb$ decomposed as the union of a $4$-ball and a regular neighborhood of the union of $S_-$ and $F$. We can also express $\CP\#\,\CPb$ as $A_{(1)}\cup M_1$ where $M_1$ is the manifold of Figure~\ref{M1N1}(a).

We pass to $A_{(1)}\cup N_1\cong \bbF_1(1)$ and 
$(T^2\x S^2)\#\, \CPb \cong \bbF_1(1)\#\,\CPb = \hA\cup N_1$. The manifold $N_1$ of Figure~\ref{M1N1}(b) is a regular neighborhood of a pair of transversely intersecting tori of self-intersections $-1$ and $0$ as described above. 
The standard surgeries on $B_T\C A\C \hA$ replace $A$ with $\TO$ and so
$(T^2\x S^2)\#\,\CPb$ is replaced with $T^4\#\,\CPb$, and these surgeries replace $\hA$ with the complement of transversely intersecting tori $\wh{T}_T$ and $T_S$ where $\wh{T}_T$
represents the homology class $[T_T]-[E]$ in $T^4\#\,\CPb$, where $E$ is the exceptional curve.

\begin{prop}\label{P3P} The standard surgeries on the three Bing double links $B_T$ in the three copies of $A\C\hA$ in $\CP\#\,3\,\CPb$ result in a symplectic $4$-manifold $Q_3$, and the core tori of the surgeries are Lagrangian in $Q_3$.\end{prop}
\begin{proof} It follows from our discussion that the under the standard surgeries, the Bing tori in $B_T\C A\C\hA$ become Lagrangian tori  in  $T^4\#\,\CPb\- (\wh{T}_T \cup T_S)$. This process carries the pinwheel structure on $\CPC$ to a pinwheel structure on a manifold $Q_3$ whose components are copies of the complements of the symplectic tori  $\wh{T}_T$ and $T_S$ in $T^4\#\,\CPb$.
It now follows from \cite{S} that $Q_3$ is a symplectic manifold with a symplectic structure extending those of the pinwheel components.
\end{proof} 
\noindent Presumably $Q_3=Sym^2(\Sig_3)$; however this has not yet been shown.

\smallskip

We next consider $\CP\#\,2\,\CPb$. As above, we write $T^4$ as the union $(\TO) \cup N$, where $N$ is a regular neighborhood of the plumbing of two disk bundles with $c_1=0$ over tori. In this case, the tori are the standard symplectic tori $T_T$ and $T_S$ where $T^4=T_T\x T_S$.  Let $\Tt$ and $\Ts$ be the copies of $T_0$ in $T_T$ and $T_S$ where we are removing disks containing the intersection point of $T_T$ and $T_S$.
Now suppose that, instead of $N$, we wish to remove from $T^4$ the complement of  $T_T$ and two parallel copies of $T_S$. This can be achieved by adjoining to $\bd(\TOst)=\bd N$ the manifold $P\x \Ts$, where $P$ is a pair of pants (Figure~\ref{Pants}). See also Figure~\ref{O} below. We have $T^4\- (T_T\cup (T_S\amalg T_S)) = (\TOst)\cup (P\x\Ts)$, and the intersection identifies $\bd \Tt\x \Ts$ with $\z\x \Ts$, where $\z$ is a boundary circle of $P$ as shown in Figure~\ref{Pants}.

\begin{figure}[ht]\vspace*{-.1in}
\begin{center}\includegraphics[scale=.85]{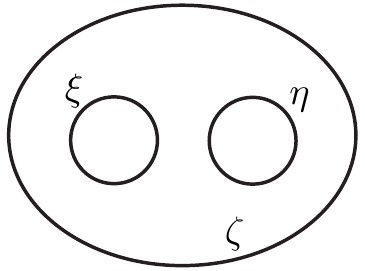}
\end{center}\vspace*{-.1in}
\caption{$P$}
\label{Pants}
\end{figure}

Instead of two parallel copies of $T_S$, we actually wish to remove a torus $2T_S$ representing the homology class $2\,[T_S]$. Then, instead of $P\x\Ts$ we obtain a nontrivial $P$-bundle, $\O$, over $\Ts$. This bundle is the restriction over $\Ts$ of the $P$-bundle over $T_S=S^1\x S^1$ whose monodromy over the first $S^1$ is the identity, and whose monodromy over the second $S^1$ is a diffeomorphism $P\to P$ which keeps $\z$ pointwise fixed and which interchanges the other two boundary circles, $\xi$ and $\eta$. In other words, $\O$ restricted over the second $S^1$-factor is $S^1$ times the complement in a solid torus of a loop going twice around the core. Then 
\[ T^4=(\TOst)\cup N=((\TOst)\cup_{\bd\Tt\x \Ts=\z\x\Ts} \O)\cup W\]
where $W$ is obtained by plumbing  two trivial disk bundles over $T^2$ twice, with both intersection points positive.

In $\TOst$ perform the standard surgeries on the tori in $\wB_T$ whose result is $A$. We obtain $A\cup (\O\cup W) = A\cup N \cong \TS$, and if we replace $N$ by $M$ we obtain $A\cup M=\SS$.

Back in $T^4$, the surfaces $T_T$ and $2T_S$ intersect at two points.
Blow up to separate the surfaces at one of these intersection points. In $N\#\,\CPb$ we have a regular neighborhood $L$ of the union of tori $T_{I,T}$ representing the homology class $[T_T]-[E]$ and $T_{I,S}$ representing $2[T_S]-[E]$ where $E$ is the exceptional curve. The neighborhood $L$ is the plumbing of two disk bundles over $T^2$ with $c_1=-1$. Write $N\#\,\CPb=R\cup L$, where $R\cong N\#\,\CPb \- L$. Then the standard surgeries on the tori in $\wB_T\C T^4\#\,\CPb$ give 
\[ \TS\#\,\CPb = A\cup (N\#\,\CPb) = A\cup (R\cup L) = I_0\cup L \]
where $I_0= A\cup R$ is the complement of a copy of $L$ in $\TS\#\,\CPb$. Thus $I_0$ is the complement in $\TS\#\,\CPb$ of tori $T'_{I,T}$ and $T'_{I,S}$ representing the homology classes $[T] - [E]$ and $2[S_T] - [E]$.  

Since $\wB_T\C \TOst$ is disjoint from $N$, the tori $T'_{I,T}$ and $T'_{I,S}$ in $\TS\#\,\CPb$ are constructed analogously to $T_{I,T}$ and $T_{I,S}$: The regular neighborhood of $T\cup S_T$ in $\TS$ may be identified with $N$.  We have $\TS=A\cup N=A\cup (\O\cup W)$, where $\O$ is now viewed as a $P$-bundle over $S_{T,0} =S_T\- D^2$.

\begin{figure}[!ht]
\begin{center}
     \subfigure[ $N$]{\includegraphics[scale=1.25]{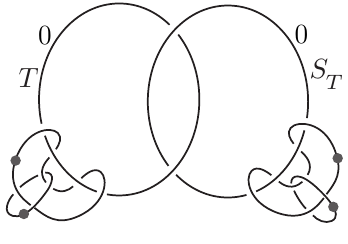}}\vspace*{-.05in}
     \hspace{.4in} 
    \subfigure[$A$]{\includegraphics[scale=1.25]{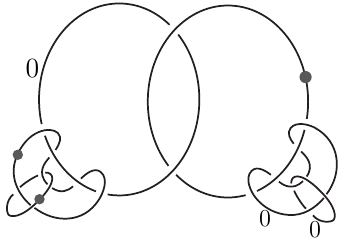}}
      \end{center}
\vspace*{-.1in}     \caption{$\TS=N\cup A$}
  \label{NuA}
\end{figure}

\bigskip

\begin{figure}[!ht]
\begin{center}
     \subfigure[ $W$]{\includegraphics[scale=1.9]{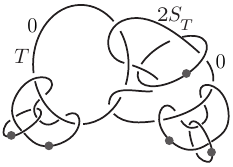}}
     \hspace{.4in} 
    \subfigure[$A\cup\O$]{\includegraphics[scale=1.9]{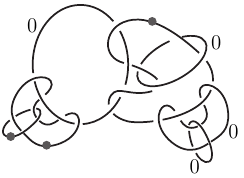}}
      \end{center}
\vspace*{-.1in}     \caption{$\TS=W\cup(A\cup\O)$}
  \label{WuAuO}
\end{figure}

\begin{figure}[ht]\vspace*{-.1in}
\begin{center}\includegraphics[scale=1.75]{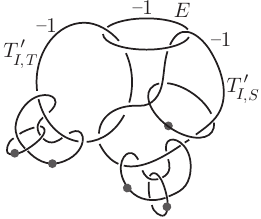}
\end{center}\vspace*{-.1in}
\caption{$W$ blown up}
\label{WHat}
\end{figure}

\begin{figure}[!ht]
\begin{center}
     \subfigure[ $L$]{\includegraphics[scale=1.3]{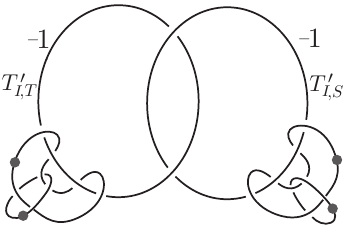}}
     \hspace{.4in} 
    \subfigure[$I_0$]{\includegraphics[scale=1.7]{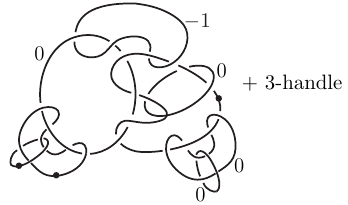}}
      \end{center}
\vspace*{-.1in}     \caption{$\TS\#\,\CPb=L\cup I_0$}
  \label{LuI}
\end{figure}

We next wish to see that this situation arises by attaching and subtracting $2$-handles starting with the complement of a pair of surfaces in $\SS\#\,\CPb$. To see this, recall the Kirby calculus depiction of $\TS = N\cup A$ as shown in Figure~\ref{NuA}.
Similarly, we have the decomposition of $\TS$ into $W\cup(A\cup\O)$ in Figure~\ref{WuAuO}. 
\noindent Next blow up $W$. This is shown in Figure~\ref{WHat}. Then Figure~\ref{LuI} shows $\TS\#\,\CPb$ as the union of $L$ and $I_0$. To accomplish this, we move the exceptional curve and the $1$-handle over to the other side of the picture. We need to make two points here. First, only part of the exceptional curve is seen in Figure~\ref{LuI}(b) because $E$ intersects both $T'_{I,T}$ and $T'_{I,S}$. Second, when the $1$-handle from Figure~\ref{WHat} is moved to the other side, it becomes a $3$-handle. 

\begin{figure}[!ht]
\centering
     \subfigure[]{\includegraphics[scale=1.9]{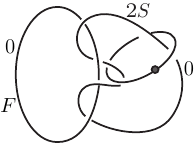}}\hspace{.25in}
    \subfigure[]{\includegraphics[scale=1.9]{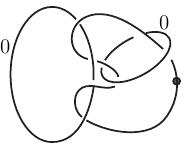}}\\
       \hspace*{.2in} \subfigure[]{\includegraphics[scale=1.25]{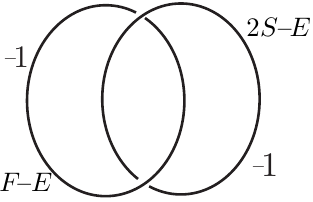}} \hspace{.25in}
    \subfigure[]{\includegraphics[scale=1.5]{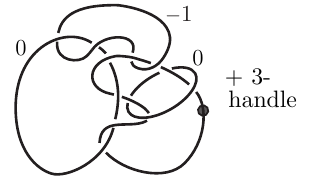}}
\vspace*{-.1in}     \caption{$(a)\cup (b)=\SS$,\ \ $(c)\cup (d)=\SS\#\,\CPb$}
  \label{SS1}
\end{figure}

By removing $1$- and $2$-handle pairs from Figures~\ref{WuAuO} and \ref{LuI}, we obtain Figure~\ref{SS1}. The union of (a) and (b) is $\SS$, and the union of (c) and (d) is $\SS\#\,\CPb$.

\begin{figure}[ht]
\begin{center}
\includegraphics[scale=1.5]{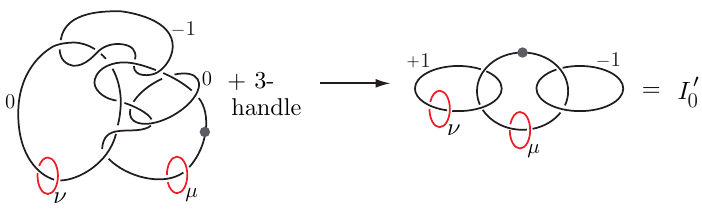}\vspace*{-.1in}
\end{center}
\caption{$I_0'$}
\label{I'}
\end{figure}

The complement of $F-E$ and $2S-E$ in $\SS\#\,\CPb$ is given in Figure~\ref{SS1}(d). In that figure, slide the left-hand $0$ over the $-1$. This is shown in Figure~\ref{I'}. (After this handle slide, the $0$-framed $2$-handle cancels with the $3$-handle.)
Figure~\ref{I'} also illustrates meridians $\mu$ and $\nu$ to the $1$- and $2$-handle. As indicated in the figure, we call this manifold $I'_0$.

\smallskip

\begin{figure}[!ht]
\begin{center}
     \subfigure[]{\includegraphics[scale=1.3]{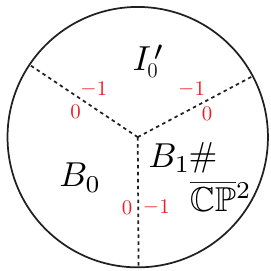}}
     \hspace{.5in} 
    \subfigure[]{\includegraphics[scale=1.3]{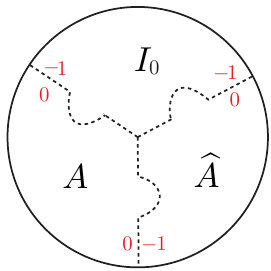}}
      \end{center}
\vspace*{-.1in} \caption{Pinwheel structures for $\CP\#\,2\,\CPb$}
  \label{PPre}
\end{figure}

\begin{figure}[ht]
\begin{center}
\includegraphics[scale=1.6]{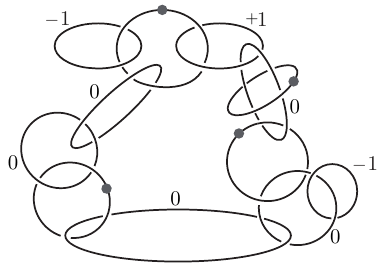}\vspace*{-.1in}
\end{center}
\caption{}
\label{start}
\end{figure}

\begin{lem}\label{P2Ppre} The $3$-fold pinwheel described by Figure~\ref{PPre}(a) is a pinwheel structure for $\SS\#\,\CPb\cong\cpc$.
\end{lem}
\begin{proof} This is an exercise in Kirby calculus. The union of the three manifolds $I_0'$, $B_1\#\,\CPb$, and $B_0$ glued in a cyclical fashion as indicated in Figure~\ref{PPre}(a) is shown in Figure~\ref{start}. What must be seen is that the union of this with $T^2\x D^2$ is $\cpc$. Handle slides, cancellations, and two blowdowns plus one anti-blowdown applied to Figure~\ref{start} reduce it to the Borromean link with two `$0$'s and a `dot'. A handle picture for $T^2\x D^2$ is the Borromean link with two `dots' and a `$0$'. In the union, all these handles cancel and we obtain $S^4$. Taking into account the blowdowns and anti-blowdown, we see that we get $\cpc$.
\end{proof}

As is indicated by the new $0$-framed $2$-handles in Figure~\ref{start}, the Handle Trading Condition is satisfied for the pinwheel structure of Lemma~\ref{P2Ppre}(a), and so Lemma~\ref{trading} gives us the pinwheel structure of Figure~\ref{PPre}(b) on $\cpc$.

\begin{prop}\label{P2P} The standard surgeries on the three Bing double links $B_T\C A\C\CP\#\,2\,\CPb$ in the pinwheel structure given by  Figure~\ref{PPre}(b)  result in a symplectic $4$-manifold $Q_2$ and the core tori of the surgeries are Lagrangian in $Q_2$.\end{prop}
\begin{proof} As in Proposition~\ref{P3P}, this follows from \cite{S}. The manifold $Q_2$ has a symplectic pinwheel structure as the $3$-fold sum of $(\tilde{X_i}; \tilde{S_i}, \tilde{T_i})$ where 
\begin{eqnarray*}
(\tilde{X_0}; \tilde{S_0}, \tilde{T_0})&=&(T^4\#\,\CPb; T_{I,T}, T_{I,S})\\
(\tilde{X_1}; \tilde{S_1}, \tilde{T_1})&=&(T^4\#\,\CPb; \wh{T}_T, T_S)\\
(\tilde{X_2}; \tilde{S_2}, \tilde{T_2})&=&(T^4; T_T, T_S) 
\end{eqnarray*}
\end{proof}

\section{Component pieces and their fundamental groups}\label{Pieces}

In this section we shall examine the fundamental groups of the the complements of our pair $B_T$ of Bing tori in the pinwheel components $A$, $\hA$, and $I_0$ of Section~\ref{PS}. One important goal is to see that when we choose our basepoint to lie on the central torus $T_c$, appropriate elements of these fundamental groups are represented by loops lying on the boundaries of these components. This will assist in the computations in Section~\ref{Exotic} where we piece together the fundamental group of a manifold from the components of a pinwheel structure on it.

\begin{figure}[ht]
\begin{center}
\includegraphics[scale=3.25]{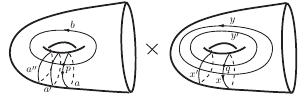}\vspace*{-.1in}
\end{center}
\caption{$\Tt\x \Ts$}
\label{T0xT0}
\end{figure}

We are helped greatly by the analysis of Baldridge and Kirk in \cite{BK}, where they study the fundamental group of the complement of the Lagrangian tori in $T_0\x T_0$ which are exactly the core tori $\L_1$ and $\L_2$ of the result of surgeries on $B_T\C A$.  Following the notation in that paper (with only minor changes) and referring to Figure~\ref{T0xT0}, the disjoint Lagrangian tori are $\L_1= a'\x x'$ and $\L_2=a''\x y'$. The relevant classes in $\pi_1$ are 
$x= \{p\}\x x$, $y= \{p\}\x y$, $a=a\x \{q\}$, and $b=b\x \{q\}$. Since $\L_1$ and $\L_2$ are Lagrangian, we refer to pushoffs of loops in $\L_i$ to parallel loops in $\L_i$ as `Lagrangian pushoffs'.

\begin{lem}[\cite{BK}]\label{bk} There are basepaths in $\Tt\x \Ts\- (\L_1\cup \L_2)=\Tt\x \Ts\- \wB_T$ to the boundaries of tubular neighborhoods of $\L_1$ and $\L_2$ so that the generators of the fundamental group $\pi_1(\Tt\x \Ts\- \wB_T, (p,q))$ that are represented by Lagrangian pushoffs of $x'$ and $a'$ on $\L_1$ are $x$ and $a$ and the meridian to $\L_1$ represents the commutator $[b^{-1},y^{-1}]$. Similarly, Lagrangian pushoffs of the loops $y'$ and $a''$ on $\L_2$ represent the elements $y$ and $bab^{-1}$, and the meridian to $\L_2$ represents $[x^{-1},b]$. Furthermore, 
$\pi_1(\Tt\x \Ts\- \wB_T, (p,q))$ is generated by the elements $x$, $y$, $a$, and $b$, and the relations $[x,a]=[y,a]=[y, bab^{-1}] =1$ (among others) hold.
\end{lem}

\begin{figure}[ht]
\begin{center}
\includegraphics[scale=.45]{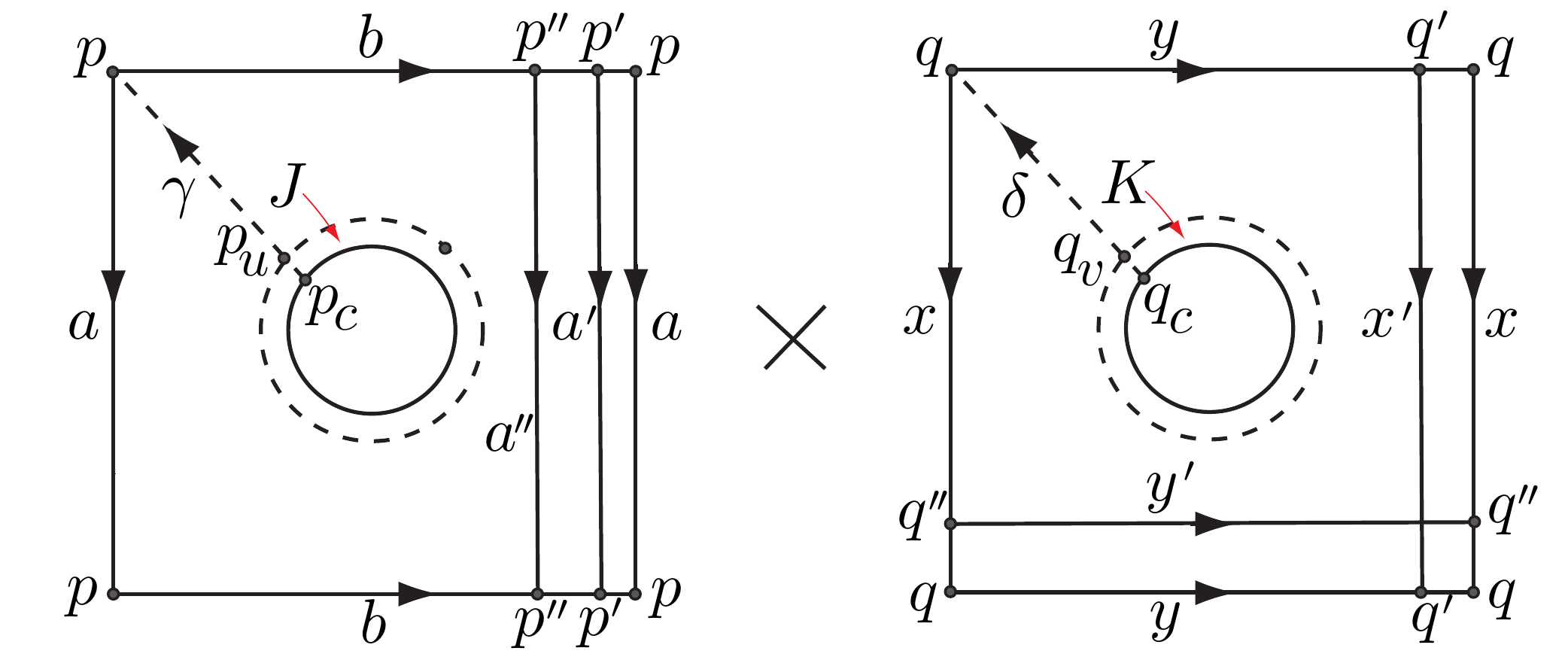}
\end{center}
\caption{$\Tt\x \Ts$}
\label{T0xT0alt}
\end{figure}

It may help in visualizing this lemma to consider (as done in \cite{RevEng}) Figure~\ref{T0xT0alt} where the punctured torus $T_0$ is viewed as a punctured disk with identifications. The loop $\{ p\}\x x'$ is based at $(p,q)$ via a path that travels backwards along $y$ from $q$ to $q'$ ({\it i.e.} $\{ p\}$ times this) and after traversing $x'$ travels back to $q$ along the same path. We base the other paths similarly. Then the lemma holds as stated.
We wish to move the basepoint to
the point $z_c=(p_c,q_c)$ which lies on the torus $\bd\, \Tt\x \bd\, \Ts$. This is done by using the path $\g\x\d$ to join $(p_c,q_c)$ to $(p,q)$. Changing the base point in this way does not affect Lemma~\ref{bk}. Since the loop $\g\ast\g^{-1}$ is nullhomotopic, we may view the loop representing $x$ based at $z_c$ as $\{p_c\}\x \xi$ where $\xi = \d\ast x\ast\d^{-1}$, and similarly for $y$, $a$, and $b$. Thus the generators of $\pi_1(\Tt\x \Ts\- \wB_T, z_c)$ are all represented by loops in $\bd(\TOst)$.

\smallskip

One obtains $A$ from $T_0\x T_0$ by surgeries on $\L_1$ (killing $x$) and $\L_2$ (killing $y$); so the lemma implies that $\pi_1(A,z_c)$ is generated by $a$ and $b$. In fact, it quickly follows from Figure~\ref{A}(a) that $\pi_1(A,z_c)$ is the free group on two generators ($a$ and $b$) which are represented by the meridians to the $1$-handles. (By `meridian' of a $1$-handle, we mean a loop parallel to its core. Thinking of a $1$-handle as a `scooped out' $2$-handle gives more intuitive meaning to this terminology.)
Notice that the classes of $a$ and $b$ generate $H_1(\Tt)$ and those of $x,y$ generate $H_1(\Ts)$.

\begin{lem} The fundamental group $\pi_1(A,z_c)$ is the free group generated by $a$ and $b$. Furthermore, $\pi_1(A\- B_T,z_c) = \pi_1(\Tt\x \Ts\- \wB_T,z_c)$, and its generators $x$, $y$, $a$, and $b$ are represented by loops on $\bd A =\bd(\TOst)$.\qed
\end{lem}

Recalling the terminology of Section~\ref{Pinwheel} (and thinking of $A$ or $\TOst$ as the $i$th component of an appropriate pinwheel)
\begin{multline*} 
\bd A=\bd(\TOst) = (\Tt\x\bd \Ts)\cup (\bd \Tt\x \Ts)=\\ ((T_i\- D^2)\x S^1) \cup T^2\x I\cup ((S_i\- D^2) \x S^1)= (T_i'\x S^1) \cup (T^2\x I) \cup (S_i'\x S^1)
\end{multline*}

In order to make the following discussion a bit easier, consider in Figure~\ref{T0xT0alt} the annuli $J=\bd\Tt\x [u,c]$ and $K=\bd\Ts\x [v,c]$. These are collars on the boundaries of the two copies of $T_0$. Let $J_c$ and $K_c$ be the boundary circles of $J$ and $K$ containing the points $p_c$ and $q_c$, and let $J_u$ and $K_v$ be the boundary circles containing $p_u$ and $q_v$.
Then $S_i'\x S^1$, the normal circle bundle over $S_i'$, is represented in this figure by 
$J_u\x (\Ts\- K)$ and $T_i'\x S^1$ is represented by $(\Tt\- J)\x K_v$.
In the above formula for $\bd A$, the 
$T^2\x I$ summand is 
$(\bd T_i'\x S^1\x [v_i,c]) \cup (\bd S_i'\x S^1\x [c,u_i])$. In the figure, this is 
$(J_c\x K) \cup ( J\x K_c)$.
In $\pi_1(\TOst\- \wB_T,z_c)$, the loop $J_c \x \{ q_c\} = \bd T_i' \x \{ q_c\}$ (in the given trivialization) represents the commutator $[a,b]$, and similarly $\{ p_c\}\x K_c = \{ p_c\}\x \bd S_i'$ represents $[x,y]$. 

We see from Figures~\ref{AA}(b) and (c) that $\hA$ is obtained from $A$ by attaching a $2$-handle with framing $-1$, and Figure~\ref{AA}(b) shows that the attaching circle for this $2$-handle is the boundary of a normal disk to the torus $T$ in $T^2\x S^2$. This is the boundary of the punctured torus $S_{T,0}$. Equivalently, after the standard surgeries on $B_T$, the attaching circle of the $2$-handle is the boundary of $\Ts$;
so in terms of our generators for $\pi_1(A\- B_T,z_c) = \pi_1(\Ts\x \Tt\- \wB_T,z_c)$, attaching the $2$-handle given by the exceptional curve adds the relation $[x,y]=1$. 

\begin{lem}\label{commutehA} The fundamental group $\pi_1(\hA,z_c)=\pi_1(A,z_c)$ and $\pi_1(\hA\- B_T,z_c)=\pi_1((\TOst)\#\,\CPb \- \wB_T,z_c)$ is obtained from $\pi_1(A\- B_T,z_c)$ by adding the relation $[x,y]=1$. Furthermore, the generators  $x$, $y$, $a$, and $b$ of the fundamental group $\pi_1(\hA\- B_T,z_c)=\pi_1( (\TOst)\#\,\CPb \- \wB_T,z_c)$ are represented by loops on $\bd \hA =\bd((\TOst)\#\,\CPb)$.\qed
\end{lem}

We next discuss the fundamental group of $I_0\- B_T$. Recall that one can construct $I_0$ by doing surgery on $\wB_T$ in $\TOst\C T^4\#\,\CPb\- (T_{I,T}\cup T_{I,S})$. 
Using the notation of the last section, $T^4\#\,\CPb\- L= (\TOst)\cup R$ where $R=N\#\,\CPb\- L$ (and the neighborhood $L$ of $T_{I,T}\cup T_{I,S}$ is the plumbing of two disk bundles over $T^2$ with $c_1=-1$). So $\pi_1(I_0\- B_T)=\pi_1((\TOst\- \wB_T)\cup R)$. We will discuss basepoints momentarily; we first need to describe $R$. To construct $R$, we begin with the manifold $\O$, the nontrivial $P$ (pair of pants)-bundle over $\Ts$ which was constructed in the last section. Referring to Figure~\ref{Pants} and below, 
\[ T^4=((\TOst)\cup_{\bd\Tt\x \Ts=\z\x\Ts} \O)\cup W\]
(Recall that $W$ is a regular neighborhood of $T_T\cup 2T_S$.)

\begin{figure}[ht]
\begin{center}
\includegraphics[scale=1]{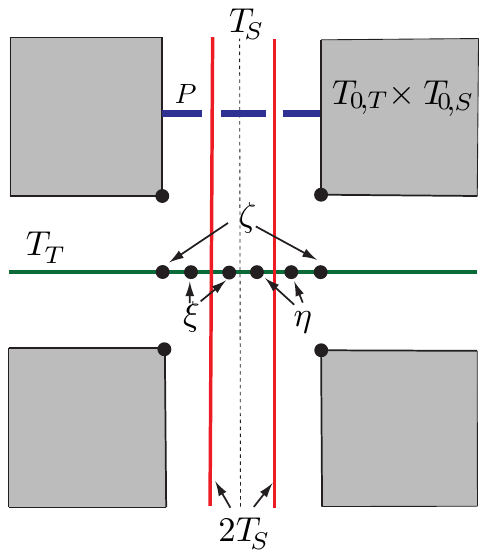}\vspace*{-.1in}
\end{center}
\caption{}
\label{O}
\end{figure}

A schematic is shown in Figure~\ref{O}, where we see $T^4=N\cup (\TOst)$. In $N$, the manifold $\O$ is the complement of a neighborhood of $T_T\cup 2T_S$. The indicated points on $\bd(\TOst)$ designate $\bd \Tt\x \bd\Ts =\z\x\bd\Ts$, the torus corresponding to the intersection point of $T_T$ and $T_S$. 
There are also tori $T_\xi = \xi\x\bd \Ts$ and $T_\eta=\eta\x\bd \Ts$ corresponding to the two intersections of the tori $T_T$ and $2T_S$.

\begin{figure}[!ht]
\begin{center}
     \subfigure[]{\includegraphics[scale=.45]{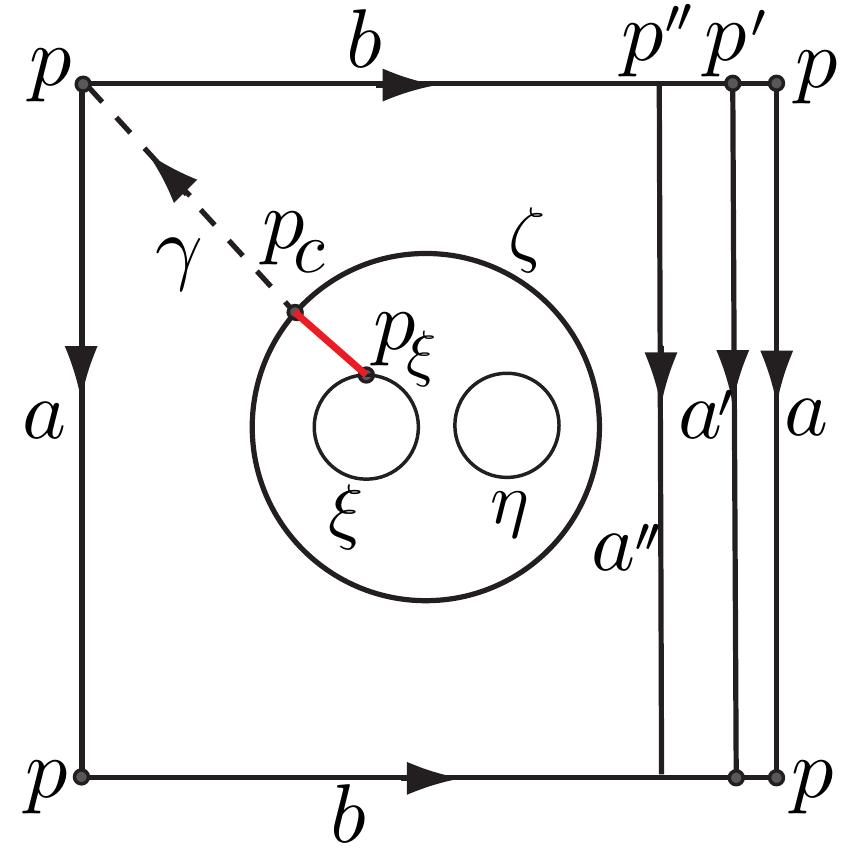}}
     \hspace{.5in} 
    \subfigure[]{\includegraphics[scale=2]{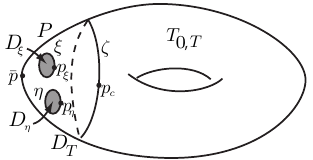}}
      \end{center}
\vspace*{-.1in}     \caption{}
  \label{T0+P}
\end{figure}

The basepoint $z_c=(p_c,q_c)$ for $\TOst$ lies on $\z\x \bd \Ts$. (See Figure~\ref{T0xT0alt}.) Use an arc from $p_c$ to $p_\xi$ inside $P\C T_T$
to move this basepoint to the point $\zx=(p_\xi,q_c)$ on $T_\xi$ for the purpose of calculating $\pi_1(T^4\- W)$ and $\pi_1(T^4\- (W\cup\wB_T))$. See Figure~\ref{T0+P}(a).

Let $x$ be the generator of $\pi_1(T_S,q_c)$ over which the monodromy of $\O$ is trivial and $y$ the generator over which the monodromy exchanges $\xi$ and $\eta$. Then 
since we have $T^4\- W = (\TOst)\cup \O$ and $z_c\in \z\x T_{0,S}=(\TOst)\cap \O$,
the groups  $\pi_1(T^4\- W,z_c)$ and $\pi_1(T^4\- (W\cup\wB_T),z_c)$ are obtained from $\pi_1(\TOst,z_c)$ and $\pi_1((\TOst)\- \wB_T,z_c)$ by adding generators $\xi$ and $\eta$, and they satisfy relations $[x,\xi]=[x,\eta]=1$, $y\xi y^{-1}=\eta$, $y\eta y^{-1}=\xi$, and $[a,b]=
\xi^{-1}\eta^{-1}$. This last relation holds because $\xi^{-1}\eta^{-1}=\z=\bd (\Tt\x \{pt\} )$. Now we can move the basepoint from $z_c$ to $\zx$ using the basepath shown in Figure~\ref{T0+P}(a); so these relations also hold in 
$\pi_1(T^4\- W,\zx)$ and $\pi_1(T^4\- (W\cup\wB_T),\zx)$

\smallskip

Let $(\bar{p},\bar{q})$ denote the intersection point $T_T\cap T_S$. Then $\Ts$ is obtained from $T_S$ by removing a disk $D_S$ containing $\bar{q}$ and similarly for $\Tt$ and a disk $D_T$ containing $\bar{p}$.
We may assume that $T_T\cap 2T_S=\{(\bar{p}_\xi,\bar{q}),(\bar{p}_\eta,\bar{q})\}$ and that $P$ is $D_T$ minus disk neighborhoods $D_\xi$ of $\bar{p}_\xi$ and $D_\eta$ of 
$\bar{p}_\eta$.

We next need to blow up to remove the intersection point $(\bar{p}_\eta,\bar{q})$ of $T_T$ and $2T_S$. For any transverse intersection of surfaces in a $4$-manifold, the result of blowing up at the intersection point is to replace a neighborhood of the intersection point, a $4$-ball containing a pair of transverse $2$-disks, with the disk bundle over $S^2$ with $c_1=-1$ containing a pair of disk fibers. In terms of the complements of the surfaces, this adds a copy of $S^1\x I\x D^2$ (the disk bundle minus two fibers) to the previous complement.

In our situation, we have $N\#\,\CPb = R\cup L$ where $L$ is a neighborhood of the tori $T_{I,T}$ (representing the homology class $[T_T]-[E]$) and $T_{I,S}$ (representing $[2T_S]-[E]$). The discussion in the above paragraph shows that $R = \O\cup (S^1\x I\x D^2)$.  Let $\eta \x I$ denote a collar in $P$ of the boundary component $\eta$.  The gluing of $S^1\x I\x D^2$ to $\O$ is given by 
\[ S^1\x I\x\bd D^2 \to   \eta\x I\x \bd\Ts, \ \ \ 
(t,r,s)\mapsto (ts^{-1},r,s)\]
I.e. $\bd D^2$ is identified with a fiber of the Hopf $S^1$-bundle, and $S^1$ is sent to $\eta\x\{q_c\}$. 

We need to apply Van Kampen's Theorem to $T^4\#\,\CPb\- L = (T^4\- W) \cup (S^1\x I\x D^2)$. Formally, the basepoint that we use needs to lie on the intersection $ \eta\x I\x \bd\Ts$; so for the purposes of this calculation we move the basepoint to $(p_\eta,q_c)$ by means of a path in $P\x \{q_c\}$. When we are done, we move the basepoint back to $\zx$ along the same path. So, in effect, we are calculating  $\pi_1(T^4\#\,\CPb\- L,\zx)$. Thus the gluing formula above implies that in $\pi_1(T^4\#\,\CPb\- L,\zx)$ and $\pi_1(T^4\#\,\CPb\- (L\cup\wB_T),\zx)$ we have $[x,y]\,\eta^{-1}=1$ and $[S^1]=\eta$.  
Thus $\eta$ lies in the image of $\pi_1(\TOst\- \wB_T,\zx)$, and since $\xi=y\eta y^{-1}$, $\xi$ also lies in the image of $\pi_1(\TOst\- \wB_T,\zx)$.
(Here we use the homotopy equivalence of $\TOst$ with $\TOst\cup (\k\x\{q_c\})$ where $\k$ is the arc from $p_c$ to $p_{\xi}$ shown in Figure~\ref{T0+P}.)

\begin{lem}\label{pi1I} The inclusion induced map: \[\pi_1(\TOst\- \wB_T,\zx)\to \pi_1(T^4\#\,\CPb\- (L\cup \wB_T),\zx)\] (equivalently, $\pi_1(A\- B_T,\zx)\to \pi_1(I_0\- B_T,\zx)$) is surjective.\qed\end{lem}

Our analysis of the blowup also shows that it changes $\bd(T^4\- W)$ to 
$\bd(T^4\#\,\CPb\- L)$ by removing $\eta\x I\x\bd\Ts$ and adding $(D_\eta\x \bd\Ts)\cup (S^1\x D_S)$.  In fact, 
\[ \bd I_0 = \bd (T^4\#\,\CPb\- L) = (T_{I,T,0}\x S^1)\cup (T_{\xi,c}\x I)\cup (S^1\x T_{I,S,0})\]
where we can identify the punctured torus $T_{I,T,0}$ with $\Tt\cup P\cup D_\eta$ and $S^1\x T_{I,S,0}$ with
$\bd_{\xi,\eta}\O\cup (S^1\x D_S)$. Here, $\bd_{\xi,\eta}\O =\bd\O \- (\z\x \Ts)$ is the inside boundary of $\O$.
The torus $T_{\xi, c}$  lies on $\bd I_0$ and contains the basepoint $\zx$. Thus the interface regions of $I_0$ are $S'\cong S^1\x T_{I,S,0}$ and $T'\cong T_{I,T,0}\x S^1$. The punctured torus $T_{I,T,0}$ has fundamental group generated by $a$ and $b$, and the punctured torus $T_{I,S,0}$ has fundamental group generated by $x$ and $y^2$, and the $S^1$ factor in $ S^1\x T_{I,S,0}$ is generated by $\xi$.

\begin{lem}\label{commute} The relation $[\xi,a]=1$ holds in the group $\pi_1(I_0\- B_T, \zx)=\pi_1(T^4\#\,\CPb\- (L\cup \wB_T), \zx)$.
\end{lem}
\begin{proof}  Since $A\- B_T\C I_0\- B_T$, Lemma~\ref{bk} implies that $[a,x]=[a,y]=1$. We have just seen that the blowup introduces the relation $\eta=[x,y]$; so
$\xi=y\eta y^{-1}= y[x,y]y^{-1}$. Thus $a$ also commutes with $\xi$.
\end{proof}

We still need to check our claim above that appropriate elements of $\pi_1(T^4\#\,\CPb\- (L\cup \wB_T),\zx) = \pi_1(I_0\- B_T,\zx)$ have representatives on $\bd(T^4\#\,\CPb\- L)=\bd I_0$. These elements are the generators $a$, $b$, $x$, and also $y^2\xi^k$  for $k\in\Z$. Before blowing up, the boundary in question is 
\[ \bd(T^4\- W) = ((\Tt\cup P)\x\bd\Ts)   \cup   \bd_{\xi,\eta}\O      \]

 Consider first $a$ and $b$. They are represented by loops in $\Tt\x\{ q_c\}$, and they become based at $\zx=(p_\xi,q_c)$ by means of a path in $P\x \{ q_c\}$ as shown in Figure~\ref{T0+P}. Thus $a$ and $b$ (based at $\zx$) have representative loops lying in $(\Tt\cup P)\x\{ q_c\}\C (\Tt\cup P)\x\bd\Ts$. These loops are unaffected by the blowup at the point $(\bar{p}_\eta,\bar{q})$.

\begin{figure}[ht]
\begin{center}
\includegraphics[scale=1.75]{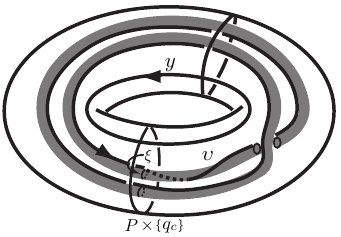}\vspace*{-.1in}
\end{center}
\caption{$\O$ restricted over $y$}
\label{Isotopey2}
\end{figure}

Next consider $x$. It was originally represented by a loop with basepoint $z_c=(p_c,q_c)$ on $\z\x\Ts$, and we used the basepath $\g\x\{ q_c\}$ to move its basepoint to $\zx$. The bundle $\O$ is trivial over $x$; so one can construct an annulus whose intersection with each fiber $P\x$(point in $x$) is the path $\g\x$(point in $x$). This gives an isotopy starting at $x$ and ending at a representative $\tilde{x}$ of $x$ based at $(p_\xi,q_c)$ which lies on $\bd_{\xi,\eta}\O$. Also, because the monodromy of the bundle $\O$ over $y$ has order two, there is a loop $\upsilon\C \bd_{\xi,\eta}\O$, depicted in Figure~\ref{Isotopey2}, which represents $y^2$. Other loops on $\bd_{\xi,\eta}\O$ which are lifts of $y^2$ are of the form $y^2\xi^k$ and are represented by $\upsilon\,\xi^k$. These loops lie on $\bd(T^4\- W)$, and our description of how blowing up changes the boundary shows that $\tilde{x}$ and $\upsilon\,\xi^k$ lie on the boundary of $T^4\#\,\CPb\- L$ after blowing up. They lie on the branch of the blown up surfaces corresponding to $\bd_{\xi,\eta}\O$; {\em{i.e.}} on
$\bd_{\xi,\eta}\O\cup (S^1\x D_S)$.

\begin{lem}\label{bdI} The elements  $a$, $b$, $x$, and $y^2\xi^k$, $k\in \Z$, of $\pi_1(T^4\#\,\CPb\- (L\cup \wB_T),\zx)=\pi_1(I_0\- B_T,\zx)$ are represented by loops on $\bd(T^4\#\,\CPb\- L) = \bd I_0$. The representatives of $a$ and $b$ are contained in  
$ T_{I,T,0}\x S^1$,  and the representatives of $x$ and $y^2\xi^k$ are contained in $S^1\x T_{I,S,0}$.\qed
\end{lem}

\section{Exotic $4$-manifolds with $b^+=1$}\label{Exotic} 

We now begin to exhibit examples which illustrate how pinwheel surgeries can be utilized to produce families of exotic $4$-manifolds. There are two advantages that this technique holds over previous methods. The first is that one is able to identify the nullhomologous tori upon which surgery is done. These are copies of $B_T\C A$ embedded in standard manifolds. The second is that there is a torus common to all the components of a pinwheel structure. Taking a basepoint in this torus can greatly simplify calculations of fundamental groups.

\subsection{$\CP\#\,3\,\CPb$} 

As a first example, we show how the construction of exotic smooth manifolds homeomorphic to $\CPC$ (cf. \cite{AP1, BK, RevEng}) fits into our framework. We have seen in Proposition~\ref{P3P} that  $\CPC$ has a $3$-fold pinwheel structure where each component is a copy of $\hA$, and that after performing the standard surgeries on the tori in the three copies of $B_T$, we obtain a symplectic manifold $Q_3$. Since each of these surgeries increases $b_1$ by one and adds a hyperbolic pair to $H_2$, we see that $Q_3$ has $b_1=6$, $b^+=7$, and $b^-=9$. We have also identified three pairs of Lagrangian tori $\wB_T$ in three copies of $T_0\x T_0\C (T^4\#\,\CPb)\- (\wh{T}_T\cup T_S)$. Each pair consists of tori $\L_1$ and $\L_2$. Recall that we obtain $\CPC$ by surgery killing the class ``$x$'' in each copy of $\L_1$ and ``$y$'' in each copy of $\L_2$. We have used the quotes here, because we really need to index the copies of $\hA$, $\L_i$, etc.

We now proceed by labelling the components of the pinwheel. The components of the pinwheel structure, which are copies of $\hA$, will be denoted $\hA_i$, the index $i$ viewed as an integer mod 3. The boundary of $\hA$ is the plumbing of two $S^1$-bundles with torus base and Euler numbers $0$ and $1$. Hence  
$\bd\hA= T_0\x S^1\cup_\vp S^1\x T_0$ for an appropriate gluing map $\vp: T^2\to T^2$. Using our previous notation (see above the Handle Trading Condition) we get for the pinwheel components that $\bd \hA_i=( S'_i\x S^1)\cup (T^2\x [-1,1]) \cup (T'_i\x S^1)$. 

Before proceeding further, we establish notation for our use of basepoints: The representative loops in $\hA_i$ for $x_i$, $y_i$, $a_i$, and $b_i$ are based at $z_{c_i}=(p_{c_i},q_{c_i})\in T_{c_i}= T^2\x{0}\C T^2\x [-1,1]$. The pinwheel components are glued together so that $T'_i\x S^1$ is identified with $S'_{i+1}\x S^1$; so we slide the basepoint for $a_i$ and $b_i$ down to $(p_{T_i},q_{T_i})\in \bd T'_i \x S^1$ by isotoping the loops off $T^2\x (-1,1)$ using the $[-1,1]$-factor. Similarly we push the basepoint for $x_i$, $y_i$ up to  $(p_{S_i},q_{S_i})\in \bd S'_i \x S^1$.

Each component $\hA_i$ is glued to $\hA_{i+1}$ identifying $T'_i\x S^1$ with $S'_{i+1}\x S^1$ so that $a_i$ is identified with $x_{i+1}$, $b_i$ is identified with $y_{i+1}$, $(p_{T_i},q_{T_i})$ is identified with $(p_{S_{i+1}},q_{S_{i+1}})$, and the corresponding basepaths are identified. This identifies the loops in $\pi_1$ that represent $a_i$ and $b_i$ with those representing $x_{i+1}$ and $y_{i+1}$. 

The union of the three copies of $T^2\x [-1,1]$ is $T^3$ and gets filled in with the manifold $T^2\x D^2$. We may suppose that all three of the points $(p_{T_i},q_{T_i})=(p_{S_{i+1}},q_{S_{i+1}})$ lie in $\{pt\}\x\bd D^2\C T^2\x D^2$, and we finally choose for our basepoint $z_0=\{pt\}\x\{0\}$ via straight lines to the center $z_0$ of $\{pt\}\x D^2$.

 In each $\hA_i$, perform surgery on the Lagrangian tori $\L_{1,i}$ and $\L_{2,i}$ as follows. On $\L_{1,i}$ perform $+1$ surgery on $x_i$ with respect to the Lagrangian framing, and on  $\L_{2,i}$ perform $+1$ surgery on $y_i$ with respect to the Lagrangian framing. It follows from \cite{ADK} that the resultant manifold will admit a symplectic structure extending the one on the complement of the tori. Furthermore, according to Lemma~\ref{bk}, these surgeries add the relations $x_i\,[b_i^{-1},y_i^{-1}]=1$ and  $y_i\,[x_i^{-1},b_i]=1$.

\begin{thm}\label{b-3} By surgeries on all  six Lagrangian tori in $\wB_{T,i}$, $i=0,1,2$ in $Q_3$, one can obtain an infinite family of mutually nondiffeomorphic smooth minimal $4$-manifolds all homeomorphic to $\CPC$. Furthermore, all these manifolds can be obtained by surgeries on six nullhomologous tori which comprise three copies of $B_T\C \CPC$.
\end{thm}
\noindent Presumably, these are the manifolds of \cite{AP1, BK, RevEng}. What is new here is the construction technique and the fact that all these manifolds are obtained by surgeries on the standard manifold $\CPC$.

\begin{proof} Let $X_1$ be the symplectic manifold which is the result of these six surgeries on $Q_3$. The fundamental group $\pi_1(X_1,z_0)$ is generated by $a_i, b_i, x_i, y_i$, $i=0,1,2$, and the relations $[x_i,y_i]=1$ (by Lemma~\ref{commutehA}), and $[x_i,a_i]=[y_i,a_i] =1$ (by Lemma~\ref{bk}) hold. Furthermore, $a_i =x_{i+1}$ and $b_i=y_{i+1}$ for $i=0,1,2$. Since $\pi_1(\hA_i\- B_{T,i})$ is generated by $a_i$, $b_i$, $x_i$, and $y_i$, we see that 
$\pi_1(X_1,z_0)$ has generators $a_i$, $b_i$ for $i=0,1,2$, and $[a_i,b_i]=1$. There is no further generator coming from gluing the components of the pinwheel because the circle $I_0\cup I_1\cup I_2$ bounds a disk in the central $T^2\x D^2$. (Each $I_i\cong [-1,1]$.)

The relations arising from the surgeries are $x_i\,[b_i^{-1},y_i^{-1}]=1$ and  $y_i\,[x_i^{-1},b_i]=1$. Translating these and the other relations, we see that $\pi_1(X_1,z_0)$ is generated by the $a_i$ and $b_i$ which satisfy
\[ [a_i,b_i]=1, \ \ [a_{i-1},a_i]=[b_{i-1},a_i] =1, \ \ a_{i-1}=[b_{i-1}^{-1},b_i^{-1}], \ \ 
b_{i-1}=[b_i,a_{i-1}^{-1}] \]
Thus $b_1= [b_2,a_1^{-1}] = [ [b_0, a_2^{-1}],a_1^{-1}]=1$, using the commutativity relations 
$[a_1,a_2]=1$ and $[b_0,a_1]=1$. Now it follows from the other relations that $\pi_1(X_1,z_0)$  is trivial.  

The infinite family of manifolds referred to in the statement of the theorem is constructed by changing the surgery on the last torus $\L_{2,2}$ so that it now kills $y_2$ times the $n$th power of the meridian 
$[x_2^{-1},b_2]$ to $\L_{2,2}$. Call this manifold $X_n$. This notation agrees with that of Theorem~\ref{SW}. Changing the surgery in this way has the effect of replacing the relation $b_1=[b_2,a_1^{-1}]$ with $b_1=[b_2,a_1^{-1}]^n$ in the calculation, which goes through just as before to show that $\pi_1(X_n)$ is trivial. Theorem~\ref{SW} now implies that the $X_n$ satisfy the conclusion of the theorem.

We have seen that there are surgeries on the three copies of $\wB_T$ in $Q_3$ that give $\CPC$ and the core tori are the three copies of $B_T$.  The composition of the inverse of these surgeries together with the surgeries on $Q_3$ that give $X_n$, give a description of surgeries on the six tori in the copies of $B_T$ in $\CPC$ that give rise to $X_n$. In other words, all of the exotic manifolds $X_n$ can be obtained by surgeries on these explicitly given nullhomologous tori in $\CPC$.

We still need to see that the manifolds $X_n$ are minimal, {\it{i.e}} that they contain no sphere of self-intersection $-1$. We first show that they have just two basic classes, $\pm k$. Since $X_n$ is homeomorphic to $\CPC$, it makes sense to talk about homology generators $h$ and $e_i$, $i=1,2,3$, where $h^2=1$ and $e_i^2=-1$. We can see from the pinwheel construction that $h-e_i$ is represented by a genus $2$ surface in $X_n$. This surface is obtained as follows: The attaching circle $\g_i$ of the $0$-framed $2$-handle in the $i^{th}$ pinwheel component for $\CPC$ is the meridian to the dotted circle giving the $1$-handle. Thus $\g_i$ is identified with the meridian to $\g_{i+1}$ in an adjacent component. Thus $\g_i$ bounds a disk in that component, and the union of this disk with the core of the $2$-handle attached to $\g_i$ is a sphere representing the homology class $h-e_i$.  
After the surgeries on the three copies of $B_T$, neither of these two $2$-disks remains. Instead, $\g_i$ bounds a punctured torus in the surgered (and handle-traded) pinwheel component $C_i$. Similarly, the meridian to $\g_{i+1}$ bounds a punctured torus in the surgered $C_{i+1}$. Thus $h-e_i$ is represented by a surface of genus $2$ in $X_n$. Similarly, each class $e_i$ is represented by a torus in $X_n$. It follows that $h = (h-e_i)+e_i$ is represented by a surface of genus $3$.

If $k$ is a basic class for (the small perturbation Seiberg-Witten invariant) on $X_n$, write $k=ah-\sum_0^2b_ie_i$. The adjunction inequality holds for the small perturbation invariant, and using the surfaces described in the paragraph above and the fact that $k^2$ must equal $(3\,\sign +2\,e)(X_n)=6$, we see that the only basic classes of $X_n$ are $\pm(3h-e_0-e_1-e_2)$. The difference of these two classes has square $24$, but if one of these manifolds failed to be minimal, it would have to have a pair of basic classes, $\k\pm e$, whose difference has square $-4$. Thus the $X_n$ are minimal.
\end{proof}

\subsection{$\CP\#\,2\,\CPb$} 

We shall next consider the construction of exotic manifolds which are homeomorphic to $\cpc$. As in Section~\ref{PS}, we can start with a pinwheel structure on $\cpc$ whose components are $I_0'$, $B_1\#\, \CPb$ and $B_0$ as in Figure~\ref{PPre}(a)
and after handle trading obtain the pinwheel structure shown in Figure~\ref{PPre}(b). The pinwheel components are glued together as follows: $C_0=I_0$ is glued to $C_1=\hA$, identifying the loop representing $a_0$ with that representing $x_1$, and $b_0$ with $y_1$. Also
$C_1=\hA$ is glued to $C_2=A$ identifying the loop representing $a_1$ with that representing $x_2$ and $b_1$ with $y_2$, and  $A$ is glued to $I_0$,
identifying $a_2$ with $x_0$ and $b_2$ with $y_0^2\xi^k$ for some $k$. This $k$ is determined by our construction in \S\ref{PS}, but its precise value will not be important to us. The basepoints in the copies of $T_c$ (recall that for $I_0$ this means the point $(p_\xi,q_c)\in T_{\xi,c}$) are all identified in the central $T^2\x D^2$ of $\cpc$ in a fashion completely analogous to the $\CPC$ case, and the basepoint for our calculation again becomes $z_0=\{pt\}\x\{0\}\in T^2\x D^2$.

We follow the pattern established above to construct examples.
After performing the standard surgeries on the tori in the three copies of $B_T$ in our pinwheel structure, we obtain the symplectic manifold $Q_2$ of Proposition~\ref{P2P} which has $b_1=6$ and the same Euler number and signature as $\cpc$. Each of the pinwheel components has fundamental group generated by the elements $a_i$, $b_i$, $x_i$ and $y_i$. 

On each $\L_{1,i}$ perform $+1$ surgery on $x_i$ with respect to the Lagrangian framing, and on  $\L_{2,i}$ perform $+1$ surgery on $y_i$ with respect to the Lagrangian framing.
Again it follows from \cite{ADK} that the resultant manifold will admit a symplectic structure extending the one on the complement of the tori, and these surgeries add the relations 
 $x_i\,[b_i^{-1},y_i^{-1}]=1$ and  $y_i\,[x_i^{-1},b_i]=1$.
We obtain a theorem analogous to Theorem~\ref{b-3}.

\begin{thm}\label{b-2} By surgeries on all six of the Lagrangian tori in $\wB_{T,i}$, $i=0,1,2$ in $Q_2$, one can obtain an infinite family of mutually nondiffeomorphic smooth minimal $4$-manifolds all homeomorphic to $\cpc$. Furthermore, all these manifolds can be obtained by surgeries on six nullhomologous tori which comprise three copies of $B_T\C \cpc$.
\end{thm}
The proof is also the same as that of Theorem~\ref{b-3} with some minor modifications.
\begin{proof} Let $X_1$ be the symplectic manifold which is the result of these six surgeries on $Q_2$. The fundamental group $\pi_1(X_1,z_0)$ is generated by $a_i, b_i, x_i, y_i$, $i=0,1,2$. After applying the identifications coming from the pinwheel structure, we see that $\pi_1(X_1,z_0)$ is generated by $a_i, b_i$, $i=0,1,2$, and $y_0$, and it satisfies the relations coming from those of $A$, $\hA$, and $I_0$ with their surgeries:
\begin{description}
\item[\hspace{.2in}$A$:]\  $a_1[b_2^{-1},b_1^{-1}]=1, \ b_1[a_1^{-1},b_2]=1, \ [a_1,a_2]=1, \ [b_1,a_2]=1$.
\item[\hspace{.2in}$\hA$:]\  $a_0[b_1^{-1},b_0^{-1}]=1, \ b_0[a_0^{-1},b_1]=1, \ [a_0,a_1]=1, \ [b_0,a_1]=1, \ [a_0,b_0]=1$. (The last equality comes from the blowup relation $[x_1,y_1]=1$.)
\item[\hspace{.2in}$I_0$:]\  $a_2[b_0^{-1},y_0^{-1}]=1, \ y_0[a_2^{-1},b_0]=1, \ [a_2,a_0]=1, \ [y_0,a_0]=1, \ 
[b_2,a_0]=1$. (The last equality holds because of Lemma~\ref{commute} and $b_2=y_0^2\xi^k$.)
\end{description}

\noindent Thus $b_0= [b_1,a_0^{-1}] = [ [b_2, a_1^{-1}],a_0^{-1}] = 1$,
using the fact that $a_0$ commutes with $b_2$ and $a_1$.
It then follows from the other relations that $\pi_1(X_1,z_0)=1$. (Recall that $\xi=y_0[x_0,y_0]y_0^{-1}$.)

The infinite family of manifolds referred to in the statement of the theorem is constructed by changing the surgery on the torus $\L_{2,2}$ so that it now kills $y_2$
 times the $n$th power of the meridian $[x_2^{-1},b_2]$ to $\L_{2,2}$. Call this manifold $X_n$. This notation agrees with that of Theorem~\ref{SW}.
Changing the surgery in this way has the effect of replacing the relation 
$b_1= [b_2,a_1^{-1}]$ with $b_1= [b_2,a_1^{-1}]^n$ in the calculation, which goes through just as before to show that $\pi_1(X_n)=1$. Theorem~\ref{SW} now implies that the $X_n$ satisfy the conclusion of the theorem.

The proof that the $X_n$ are minimal follows an argument similar to that in the proof of Theorem~\ref{b-3}. Once again we see from the pinwheel structure that $X_n$ contains genus $2$ surfaces representing the classes $h-e_1$ and $h-e_2$ and a genus one surface representing $e_1$, say; so there is again a genus $3$ surface representing $h=(h-e_1)+e_1$. We see that the only possible basic classes are $\pm(3h-e_1-e_2)$, and the proof proceeds as before.
\end{proof}

\section{More exotic rational surfaces}\label{More}

\begin{figure}[!ht]
\begin{center}
     \subfigure[ $A_{(n)} $]{\includegraphics[scale=1.3]{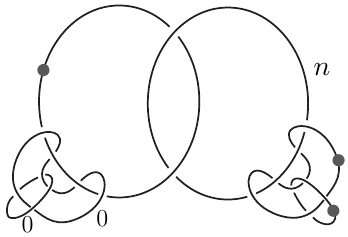}}
     \hspace{.3in} 
    \subfigure[$\hA_{(n)}$]{\includegraphics[scale=1.3]{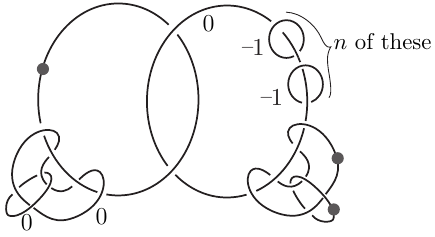}}
      \end{center}
\vspace*{-.1in}     \caption{}
  \label{AA_n}
\end{figure}

In the last section we showed how to use pinwheel surgery to obtain infinite families of manifolds homeomorphic but not diffeomorphic to $\CP\#\,k\,\CPb$ for $k=2,3$. In this section, we outline how to use similar techniques to construct families for $k=4,\dots,7$ and $k=9$, and we explain why we are having difficulty seeing  a similar construction for $k=8$. We begin by generalizing the construction of the manifold $\hA$. Recall the notation $B_n$ from Figure~\ref{Hopf}. For $n>0$, $B_n$ is naturally the complement of the negative section and a fiber in the ruled surface $\bbF_n$. 
In accordance we also consider the manifolds  $A_{(n)} $ and 
$\hA_{(n)}$ of Figure~\ref{AA_n}. In terms of this notation, our manifold $\hA =\hA_{(1)}$. Notice that $A\C \hA_{(n)}$.

\begin{figure}[!ht]
\begin{center}
     \subfigure[ Neighborhood of $S_- \cup (S_++F)$]{\includegraphics[scale=2.5]{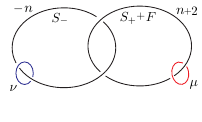}}
     \hspace{.3in} 
    \subfigure[$K_n$]{\includegraphics[scale=2.5]{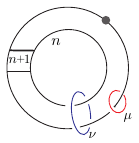}}
      \end{center}
\vspace*{-.1in}     \caption{$\bbF_n$}
  \label{Fn}
\end{figure}

We introduce next another family of manifolds which we shall need. Let $K_n$ denote the complement in $\bbF_n$ of $S_-$, the negative section, and $S_++F$, a positive section plus a fiber.  These are embedded complex curves in $\bbF_n$ which meet transversely in a single point and have self-intersections $-n$ and $n+2$. It follows that $K_n$ is the rational ball of \cite{rat} whose boundary is the lens space $L((n+1)^2,-n)$, and the fundamental group $\pi_1(K_n)=\Z_{n+1}$. (The manifold $K_0 \cong B^4$. See Figure~\ref{Fn}(b).)

Figure~\ref{Fn} displays $\bbF_n$ as the the union of a neighborhood of the intersecting spheres $S_-$ and $S_++F_n$ and $K_n$. (The `$n+1$' indicates $n+1$ full right hand twists.) The loops $\mu$ and $\nu$ are the boundaries of normal disks to the two surfaces, and they are shown in Figure~\ref{Fn}(b) lying in $\bd K_n$. It is important to note that $\mu$ and $\nu$ are isotopic in $K_n$. We can see this in Figure~\ref{Fn}(b) since there is an isotopy of $\nu$ across a belt disk of the $2$-handle. Alternatively, a fiber of $\bbF_n$ intersects both $S_-$ and $S_++F$ transversely in single points, and the complementary annulus in the fiber gives the isotopy.

We next describe an operation that will lead to a mild generalization of handle trading. Consider $K_n$ and abstractly attach a pair of $2$-handles to the Bing double of $\mu$ with both framings equal to $0$. Call this new manifold $K_n^+$. Since $\nu$ is isotopic to $\mu$ in $K_n$, the Bing double of $\nu$ is isotopic to the Bing double of $\mu$ in $K_n$, and using this isotopy and the cores of the new $2$-handles in $K_n^+$, we get disks in $K_n^+$ with boundary the Bing double of $\nu$, and we can remove neighborhoods of these disks from 
$K_n^+$ to obtain the manifold $K_n^0$. 

\begin{figure}[ht]
\begin{center}
\includegraphics[scale=1.5]{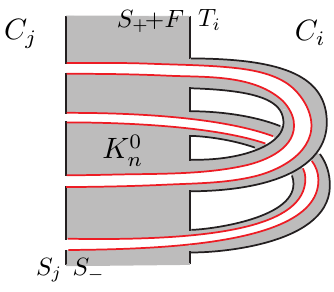}\vspace*{-.1in}
\end{center}
\caption{Pushing Through}
\label{Push}
\end{figure}

This procedure will be used as follows: Visualize $K_n$ as a pinwheel component, and consider the adjacent pinwheel component $C_i$ with interface surface $T_i$ whose meridian $m_{T_i}$ is identified with the meridian $\mu$ of $S_++F$. Suppose that the meridian $m_{T_i}$ satisfies the Handle Trading Condition; so in the pinwheel, we can add the pair of $2$-handles to $K_n$ (removing them from $C_i$)  in order to construct $K_n^+$. 
Next consider the pinwheel component $C_j$ on the other side of $K_n$ with interface surface $S_j$ whose meridian $m_{S_j}$ is identified with the meridian $\nu$ of $S_-$. We can view the construction of $K_n^0$ from $K_n^+$ as handle trading. In other words, the Handle Trading Condition is not satisfied at the $S_-$ interface of $K_n$, but after trading handles at the $S_++F$ interface, it becomes possible to trade handles in $K_n^+$. See Figure~\ref{Push}. Note that the interface surfaces of the new pinwheel component $K_n^0$ are tori (as usual after handle trading). We call this operation {\it{pushing through}}.

\begin{lem} \label{passthru} The manifold $K_n^0$ is the complement in the irrational ruled surface $\bbF_n(1)$ of symplectic tori representing the negative section, $S'_-$, and a positive section plus a fiber, $S'_++F'$.
\end{lem}
\begin{proof}  To change $\bbF_n$ into $\bbF_n(1)$, we need to remove a neighborhood of a fiber, $F\x D^2\cong S^2\x D^2$, and replace it with $F \x T_0$. We are working with the smaller manifold, $K_n$; so the change we would need to make is to replace  $S^1 \x [0,1] \x D^2$ with  $S^1 \x [0,1] \x T_0$. Thus our discussion is local to this $S^1 \x [0,1] \x D^2$ in $K_n$.

The isotopy in $K_n$ used for pushing through sits inside of a fiber of $\bbF_n$. Its trace is the annulus $S^1 \x [0,1]$. The pushing through procedure adds a pair of $2$-handles to the link $L_{\mu}=L \x {1}$ where $L$ is the Bing double of $S^1 \x {0}$ in $S^1 \x D^2$.  Let $D_1$, $D_2$ be the core disks of these $2$-handles.

We then bore out a neighborhood of $(L \x [0,1])  \cup  D_1 \cup  D_2$. The claim is that the resultant manifold is $S^1 \x [0,1] \x T_0$. Notice that each  $S^1 \x \{ t\} \x D^2$ has been replaced by the $3$-manifold obtained by adding a pair of $2$-handles to 
$L \x \{ t\}$, {\it{i.e.}} by performing $0$-surgery on both components of this link.  This is $S^1 \x T_0$ as required.
\end{proof}

As we have mentioned above, the fundamental group of $K_n$ is $\Z_{n+1}$, and it is generated by either of the meridians $\mu$, $\nu$. To construct $K_n^0$ we attach a pair of $2$-handles and then subtract a pair of $2$-handles. The handle additions add no generators to $\pi_1$ and the handle subtractions add a pair of generators $a$, $b$, which are identified with corresponding generators of the form $x_i$, $y_i$ in the adjacent pinwheel component $C_i$.

\subsection{$\CP\#\,4\,\CPb$}  
\begin{figure}[!ht]
\begin{center}
     \subfigure[]{\includegraphics[scale=1.3]{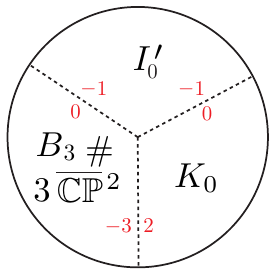}}
     \hspace{.5in} 
    \subfigure[]{\includegraphics[scale=1.3]{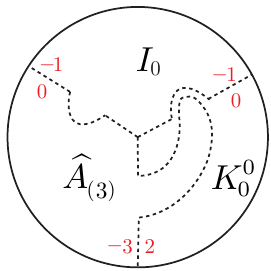}}
      \end{center}
\vspace*{-.1in} \caption{Pinwheel structures for $\CP\#\,4\,\CPb$}
  \label{PP4}
\end{figure}

The pinwheel structure we shall use to construct families for $\CP\#\,4\,\CPb$ is given in Figure~\ref{PP4}. 

\begin{figure}[ht]
\begin{center}
\includegraphics[scale=1.6]{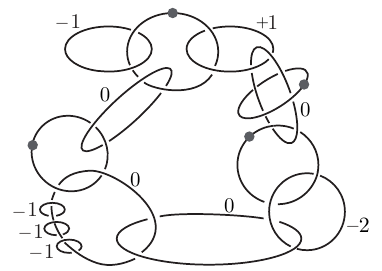}\vspace*{-.1in}
\end{center}
\caption{}
\label{start_4}
\end{figure}

The interface surface between the pinwheel components $I_0'$ and $B_3\# 3\,\CPb$ satisfies the Handle Trading Condition because the corresponding interface surface of $B_3\# 3\,\CPb$ is a fiber of $\bbF_3$ (the other interface surface of this component is a negative section), and the disk used for handle trading is provided by the positive section blown up three times. It is the core of the handle labelled $n$ ($=3$) in Figure~\ref{Hopf} after the blowups. The interface surface between the components $K_0$ and $I'_0$ also satisfies the Handle Trading Condition using the normal disk to the $2$-handle with framing $+1$ in $I'_0$. 
To see that this pinwheel structure actually gives $\CP\#\,4\,\CPb$, one argues as in Lemma~\ref{P2Ppre}. The starting Kirby calculus diagram is shown in Figure~\ref{start_4}. 

\smallskip

The model symplectic manifold $Q_4$ is built by performing standard surgeries on the copies of $B_T$ in $I_0$ and $\hA_{(3)}$. No surgeries are performed in $K_0^0$. Thus (employing Lemma~\ref{passthru}) $Q_4$ has the symplectic pinwheel structure given by 
\begin{eqnarray*} 
(\tilde{X_0}; \tilde{S_0}, \tilde{T_0})&=&(T^4\#\,\CPb; T_{I,T}, T_{I,S})\\
(\tilde{X_1}; \tilde{S_1}, \tilde{T_1})&=&(\TS; T_{S_++F}, T_{S_-})\\
(\tilde{X_2}; \tilde{S_2}, \tilde{T_2})&=&(T^4\#3\,\CPb; T_{T,3}, T_S)
\end{eqnarray*}
where $T_{T,3}$ is a symplectic torus representing $T_T-E_1-E_2-E_3$ and in $T^2\x S^2=\bbF_0(1)$, $T_{S_++F}$ and $T_{S_-}$ are tori representing the obvious classes.

To construct the infinite family of mutually nondiffeomorphic manifolds all homeomorphic to $\CP\#\,4\,\CPb$, we follow the proof of Theorem~\ref{b-2}, except that here we only surger the tori $\L_{0,i}$ and $\L_{2,i}$ coming from the pinwheel components $I_0$ and $\hA_{(3)}$. The component $K_0^0$ is diffeomorphic to $\TS\- (T_{S_++F}\cup T_{S_-})$ where $\pi_1(\TS,z_c)$ is generated by $a_1$ and $b_1$ (where we identify $K_0^0$ with the component $C_1$). So $\pi_1(K_0^0,z_c)$ is normally generated by $a_1$, $b_1$, and a meridian $\mu$ to 
$T_{S_++F}$.

\begin{thm}\label{b-4} By surgeries on the four Lagrangian tori comprising the copies of $\wB_{T,i}$, $i=0,2$ in $Q_4$ which live in $\tilde{X_0}$ and $\tilde{X_2}$ in $Q_4$, one can obtain an infinite family of mutually nondiffeomorphic minimal smooth $4$-manifolds all homeomorphic to $\CP\#\,4\,\CPb$. Furthermore, all these manifolds can be obtained by surgeries on four nullhomologous tori which comprise two copies of $B_T\C I_0\cup \hA_{(3)}\C\CP\#\,4\,\CPb$.
\end{thm}
\begin{proof} We first show that the manifold $X$ obtained by performing $+1$ surgeries on the Lagrangian tori $\L_{1,i}$ and $\L_{2,i}$ in $\tilde{X_i}$, $i=0,2$, is simply connected.  Our pinwheel components are 
$C_0=I_0$, $C_1= K_0^0$, and $C_2=\hA_{(3)}$. The fundamental group of $\hA_{(3)}\- B_T$ is the same as that of $\hA\- B_T$ since the extra blowups do not change $\pi_1$. The pinwheel construction identifies $a_2$ with $x_0$ and $b_2$ with $y_0^2\,\xi^k$. Furthermore, the pushing through process identifies $a_0$ and $b_0$ with the belt circles of the $2$-handles first attached to $K_0$ to construct $K_0^+$, and these are in turn identified with $a_1=x_2$ and $b_1=y_2$ after pushing through. Thus $a_0=a_1$ and $b_0=b_1$.

Analogous to the proof of Theorem~\ref{b-2},  we see that $\pi_1(X,z_0)$ is generated by $a_i, b_i$, $i=0,2$, $y_0$, and the meridian $\mu$ from $K_0^0$, and it satisfies the relations coming from those of $\hA_{(3)}$, and $I_0$ with their surgeries:
\begin{description}
\item[\hspace{.15in}$\hA_{(3)}$:]\  $a_0[b_2^{-1},b_0^{-1}]=1, \ b_0[a_0^{-1},b_2]=1, \ [a_0,a_2]=1, \ [b_0,a_2]=1, \ [a_0,b_0]=1$.
\item[\hspace{.15in}$I_0$:]\  $a_2[b_0^{-1},y_0^{-1}]=1, \ y_0[a_2^{-1},b_0]=1, \ [a_2,a_0]=1, \ [y_0,a_0]=1, \ 
[b_2,a_0]=1$ . (Recall that this last equality holds because of Lemma~\ref{commute} and $b_2=y_0^2\,\xi^k$.)
\end{description}

\noindent  Thus $b_0=[b_2,a_0^{-1}]=1$, and then it follows that $a_0$, $a_2$, and $y_0$ vanish, and so we also have $\xi=y_0[x_0,y_0]y_0^{-1}=y_0[a_2,y_0]y_0^{-1}=1$; hence $b_2=1$.  Finally, $\mu=\nu$ is identified with the boundary of the meridian to the $1$-handle of $\hA_{(3)}$ (see Figure~\ref{AA_n}(b)). Thus $\mu=[a_2,b_2]=1$; so $\pi_1(X,z_0)=1$.

It follows as usual that the manifolds $X_n$, obtained by performing $1/n$-surgery rather than $+1$-surgery on the last torus $\L_{2,2}$, are also simply connected. To see that these manifolds are minimal, we again use an argument based on the adjunction inequality is in the proof of Theorem~\ref{b-3}. This is easy once we see that these manifolds have embedded genus $2$ surfaces representing the classes $h-e_i$, $i=1,\dots,4$ and an embedded  genus $3$ surface representing $h$.

We shall describe these surfaces referring to Figures~\ref{PP4} and \ref{start_4}. The class $e_4$ is represented as follows. The $2$-handle of $I_0'$ which is labelled `$-1$' (on the upper left in Figure~\ref{start_4}) gives a punctured sphere, punctured by the $1$-handle of $I'_0$. The meridian of that $1$-handle is isotopic to (identified with) the meridian of the $1$-handle in $B_3\#\,3\,\CPb$. This meridian is, in turn, isotopic to the attaching circle of the $2$-handle labelled `$0$' in $B_3\#\,3\,\CPb$. Thus we see a sphere of self-intersection $-1$ in $\CP\#\,4\,\CPb$. After the surgeries on the tori of the copy of $B_T$ in $\hA_{(3)}$ this sphere is no longer extant, but a torus representing $e_4$ is, since normal disks to the $1$-handle become punctured tori. 

Instead of starting with the $2$-handle of $I_0'$ which is labelled `$-1$', we could begin with the handle labelled `$+1$'. Then this construction produces a $2$-sphere in $\CP\#\,4\,\CPb$ representing $h$. Notice that both classes $h$ and $e_4$ intersect the spheres of self-intersection $-1$ which come from the three $2$-handles labelled `$-1$' in $B_3\#\,3\,\CPb$. It follows easily that these three spheres represent the classes $h-e_i-e_4$, $i=1,2,3$. Since the surgeries on $B_T\C \hA_{(3)}$ turn normal disks to the `$0$' into punctured tori, the classes 
$h-e_i-e_4$, $i=1,2,3$ are represented by tori in the surgered manifolds $X_n$.
It follows that $h-e_i= (h-e_i-e_4)+e_4$, $i=1,2,3$, are represented by genus $2$ surfaces in $X_n$. 

Returning to $\CP\#\,4\,\CPb$, there is a sphere representing a class of self-intersection $0$ arising from the $2$-handle labelled $0$ at the bottom of Figure~\ref{start_4}. (Alternatively, this sphere is formed from the cocores of the $2$-handle labelled $0$ in $B_3\#\,3\,\CPb$ and the $2$-handle labelled $-2$ in $K_0$, since the pinwheel construction identifies their boundaries.) Since this sphere intersects both $h$ and $e_4$ once and is orthogonal to 
$h-e_i-e_4$, $i=1,2,3$, it represents $h-e_4$. Surgery plus the pushing through operation make this surface genus $2$. Finally, $h=(h-e_4)+e_4$ is now seen to be represented by a genus $3$ surface in $X_n$, and we may carry out the adjunction inequality argument as planned to see that $X_n$ is minimal.
\end{proof}

\subsection{$\CP\#\,5\,\CPb$}  
\begin{figure}[!ht]
\begin{center}
     \subfigure[]{\includegraphics[scale=1.3]{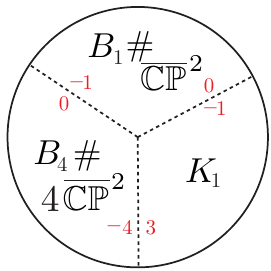}}
     \hspace{.5in} 
    \subfigure[]{\includegraphics[scale=1.3]{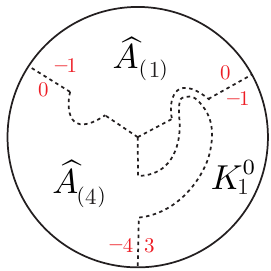}}
      \end{center}
\vspace*{-.1in} \caption{Pinwheel structures for $\CP\#\,5\,\CPb$}
  \label{PP5}
\end{figure}

The construction of families for $\CP\#\,5\,\CPb$ follows our established pattern, using the pinwheel structure of Figure~\ref{PP5}. One sees that these pinwheels give $\CP\#\,5\,\CPb$ via our usual Kirby calculus argument.

We get our infinite family of manifolds as before. The argument is parallel to, but easier than, that of Theorem~\ref{b-4}. 

\subsection{$\CP\#\,6\,\CPb$} \label{6}
\begin{figure}[!ht]
\begin{center}
     \subfigure[]{\includegraphics[scale=1.3]{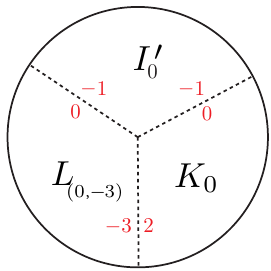}}
     \hspace{.5in} 
    \subfigure[]{\includegraphics[scale=1.3]{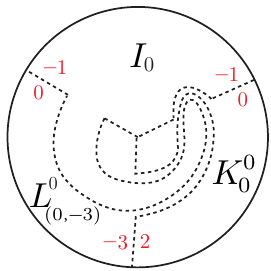}}
      \end{center}
\vspace*{-.1in} \caption{Pinwheel structures for $\CP\#\,6\,\CPb$}
  \label{PP6}
\end{figure}

This construction is similar to that for $\CP\#\,4\,\CPb$. In that construction (Figure~\ref{PP4}), replace the pinwheel component $B_3\#\,3\,\CPb$ with the manifold $L_{(0,-3)}$ which is obtained by removing the spheres $S_++F-E_1-E_2-E_3$ and $S_--E_4-E_5$ from $\bbF_1\#\, 5\,\CPb$. These spheres have self-intersections $0$ and $-3$ as required. See Figure~\ref{PP6}. Note that, as for the manifolds $K_n$, the meridians to the spheres $S_++F-E_1-E_2-E_3$ and $S_--E_4-E_5$ are isotopic in $L_{(0,-3)}$. (This isotopy is given by the fiber of $\bbF_1$.) 

In this case we need to push through two pinwheel components. Surgery is then performed on the single pair of Bing tori in $I_0$. The notation $L^0_{(0,-3)}$ is meant to be analogous to $K_n^0$ above. The proof of Lemma~\ref{passthru} shows that $L^0_{(0,-3)}$  is the complement in the blown up irrational ruled surface $\bbF_n(1)\#\,5\,\CPb$ of symplectic tori representing $S'_++F-E_1-E_2-E_3$ and $S'_--E_4-E_5$.
Standard surgeries on the tori in $B_T \C I_0$ give rise to the symplectic manifold $Q_6$ which has the symplectic pinwheel structure given by 
\begin{eqnarray*} 
(\tilde{X_0}; \tilde{S_0}, \tilde{T_0})&=&(T^4\#\,\CPb; T_{I,T}, T_{I,S})\\
(\tilde{X_1}; \tilde{S_1}, \tilde{T_1})&=& (\TS; T_{S_++F}, T_{S_-}) \\
(\tilde{X_2}; \tilde{S_2}, \tilde{T_2})&=&(\bbF_1(1)\#\,5\,\CPb; T_{S_++F-E_1-E_2-E_3}, T_{S_--E_4-E_5})
\end{eqnarray*}
with notation as above. See the $\CP\#\,7\,\CPb$ case below for comments on a similar fundamental group calculation. Notice that our construction leaves none of the exceptional curves from the blowups, and our usual adjunction inequality argument shows that we get a family of minimal $4$-manifolds.

\subsection{$\CP\#\,7\,\CPb$} 
\begin{figure}[!ht]
\begin{center}
     \subfigure[]{\includegraphics[scale=1.3]{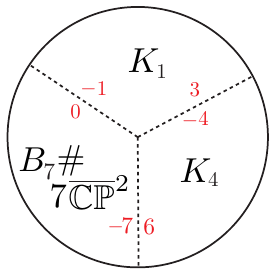}}
     \hspace{.5in} 
    \subfigure[]{\includegraphics[scale=1.3]{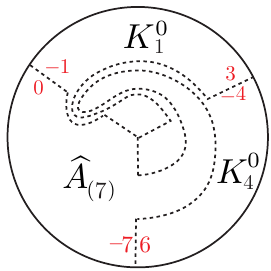}}
      \end{center}
\vspace*{-.1in} \caption{Pinwheel structures for $\CP\#\,7\,\CPb$}
  \label{PP7}
\end{figure}

Use the pinwheel structure given in Figure~\ref{PP7}. Again we need to push through two pinwheel components. Surgery is performed on a single pair of Bing tori.

Standard surgeries on the tori in $B_T \C \hA_{(7)}$ give rise to the symplectic manifold $Q_7$ which has the symplectic pinwheel structure given by 
\begin{eqnarray*} 
(\tilde{X_0}; \tilde{S_0}, \tilde{T_0})&=&(\bbF_1(1); T_{S_++F}, T_{S_-})\\
(\tilde{X_1}; \tilde{S_1}, \tilde{T_1})&=&(\bbF_4(1); T_{S_++F}, T_{S_-})\\
(\tilde{X_2}; \tilde{S_2}, \tilde{T_2})&=&(T^4\#7\,\CPb; T_{T,7}, T_S)
\end{eqnarray*}
where $T_{T,7}$ is a symplectic torus representing $T_T-\sum_{i=1}^7E_i$, and the other notation is as above.

\begin{thm}\label{b-7} By surgeries on the two Lagrangian tori comprising the copy of $\wB_T \C \tilde{X_2}$, in $Q_7$, one can obtain an infinite family of mutually nondiffeomorphic minimal smooth $4$-manifolds all homeomorphic to $\CP\#\,7\,\CPb$. Furthermore, all these manifolds can be obtained by surgeries on two nullhomologous tori which comprise two copies of $B_T\C \hA_{(7)}\C\CP\#\,7\,\CPb$.
\end{thm}
\begin{proof} We first show that the manifold $X$ obtained by performing $+1$ surgeries on the Lagrangian tori $\L_{1,2}$ and $\L_{2,2}$ in $\tilde{X_2}$ is simply connected.
The pushing through process identifies $a_2$ and $b_2$ with the belt circles of the $2$-handles first attached to $K_1$ to construct $K_1^+$, and these are in turn identified with $a_0$ and $b_0$. The second pushing through operation identifies these classes with $a_1$ and $b_1$. Thus $a_0=a_1=a_2$ and $b_0=b_1=b_2$. Therefore,
$\pi_1(X;z_0)$ is generated by $a=a_i$, $b=b_i$, and meridians $\mu_0=\nu_0$ of $S_++F$ and $S_-$ in $\bbF_1$ and $\mu_1=\nu_1$ of $S_++F$ and $S_-$ in $\bbF_4$.

The usual relations coming from $\hA_{(7)}$ in this case reduce quickly to $a=1$, $b=1$. As in the proof of Theorem~\ref{b-4} we have $\mu_i=[a,b]=1$, $i=0,1$; so $\pi_1(X,z_0)=1$, and the rest of the proof follows analogously to the examples above.
\end{proof}

\subsection{$\CP\#\,8\,\CPb$ and $\CP\#\,9\,\CPb$} These constructions are slightly different from those above. 
We begin by explaining why the type of pinwheel surgeries that we have already described will not work for these manifolds. In each of our pinwheel surgery constructions, the manifold $X_1$ obtained by performing only $+1$ surgeries is symplectic and not diffeomorphic to $\CP\#\,k\,\CPb$. The adjunction formula shows that its canonical class is represented by a symplectic surface of genus $10-k$. The canonical class of the standard manifold  $\CP\#\,k\,\CPb$ has genus $1$, and each surgery on a Bing torus pair $B_T$ in a copy of $A$ increases the genus of the disks which span meridians to the `$0$' and the `dot' by one. For example, letting $k=3$, the canonical class of $\CPC$ is given by $-3h+e_ 1+e_2+e_3$, and in each pinwheel component $B_1\#\,\CPb$, it is represented by a pair of transverse disks normal to the `$0$' and the `dot'. Gluing these together along the interfaces, they add up to three $2$-spheres representing the classes $-h+e_i$, and arranged in a cyclical fashion; so that smoothing their intersection points gives a torus representing  $-3h+e_ 1+e_2+e_3$.
The surgeries increase the genus of each of these surfaces to $2$, and the result after smoothing double points is the canonical class of $X_1$ of genus $7= 10-k$. The other examples can be analyzed in a like fashion.

We notice that each pinwheel component that is surgered must increase the genus of a symplectic representative of the canonical class by at least $2$; so the minimal genus which can be achieved by this technique is $3$. However,  the canonical class of a symplectic manifold homeomorphic but not diffeomorphic to $\CP\#\,8\,\CPb$ has genus $2$; so it cannot be obtained by this particular type of pinwheel surgery.  

\begin{figure}[!ht]
\begin{center}
     \subfigure[$\CP\#\,8\,\CPb$]{\includegraphics[scale=1.3]{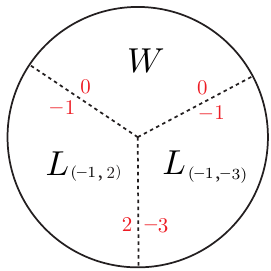}}
     \hspace{.5in} 
    \subfigure[$\CP\#\,9\,\CPb$]{\includegraphics[scale=1.3]{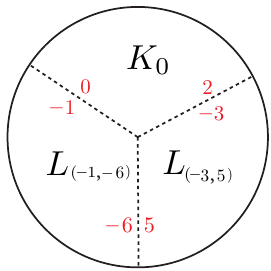}}
      \end{center}
\vspace*{-.1in} \caption{Pinwheel structures for $\CP\#\,8\,\CPb$ and $\CP\#\,9\,\CPb$}
  \label{PP89}
\end{figure}

However, one does have different pinwheel constructions for $\CP\#\,8\,\CPb$ and $E(1)=\CP\#\,9\,\CPb$ as in Figures~\ref{PP89}(a) and (b). The manifolds $L_{(j,k)}$ are obtained by modifying the constructions of the $K_n$ above. For example, $L_{(-1,-3)}$ is the complement of $S_-$ and $S_++F-E_1-\dots-E_6$ in $\bbF_1\#\,6\,\CPb$. The component $W$ is the complement of $S$ and $2S+F-2E$ in $\bbF_0\#\,\CPb$.

We leave it as a simple exercise for the reader to see representatives for the canonical classes in these pinwheel structures, tori of self-intersections $1$ and $0$, respectively. Surgery on the Bing double of this last torus was the subject of the example in Section~\ref{Bing}. By replacing $K_0$ with $K_0^0$ and each $L_{(j,k)}$ with $L_{(j,k)}^0$; {\it{i.e.}} by replacing the pinwheel components of $E(1)=\CP\#\,9\,\CPb$ with components obtained in like fashion from ruled surfaces over $T^2$ rather than over $S^2$, we obtain the model manifold $Q_9$, and it is easy to see that $Q_9$ is the elliptic surface over $T^2$ obtained by fiber summing $E(1)$ with the product elliptic surface $T^2\x T^2$. Thus $Q_9$ is obtained from $E(1)$ by removing a neighborhood of a fiber, $T^2\x D^2$ and replacing it with $T^2\x T_0$. In other words, we have performed standard surgeries on the standard copy of $B_T$ in the neighborhood of a fiber in $E(1)$ in order to obtain $Q_9$. The surgeries on the corresponding Lagrangian tori in $E(1)$ will give exactly the examples of Section~\ref{Bing}.

We have not as yet been able to effect a construction of exotic manifolds homeomorphic to $\CP\#\,8\,\CPb$ utilizing pinwheels. However in \cite{Santeria} we constructed a manifold $R_8$ homeomorphic to $\CP\#\,8\,\CPb$ and with trivial Seiberg-Witten invariant (as for $\CP\#\,8\,\CPb$) such that $R_8$ contains nullhomologous tori upon which surgery yields an infinite family of exotic manifolds homeomorphic to $\CP\#\,8\,\CPb$.

\subsection{The main theorem}\label{Main} The upshot of this section and the last is the following theorem:

\begin{thm} For $k=2,\dots,7$ and $k=9$, there are nullhomologous tori embedded in $\CP\#\,k\,\CPb$ upon which surgery gives rise to an infinite family of mutually nondiffeomorphic minimal $4$-manifolds homeomorphic to $\CP\#\,k\,\CPb$.\qed
\end{thm}

\section{Final Remarks}

The fact that surgery on a well-chosen collection of (nullhomologous) tori in a fixed smooth manifold can alter smooth structures gives rise to a target (and optimistic) classification scheme for simply connected smooth $4$-manifolds.

\begin{quest} Is it possible that every simply connected smooth $4$-manifold is obtained by surgery on tori in a connected sum of copies, with either orientation, of $S^{4}$, $\CP$, $\CPb$, $S^{2}\times S^{2}$, and the 
$K3$ surface?
\end{quest}

\subsection*{Acknowledgements}  R.F. was partially supported by NSF Grants DMS-0704091 and DMS-1006322 and R.J.S. was supported by no federal research dollars.


\begin{thebibliography}

\bibitem[1]{ABP} A. Akhmedov, R.I. Baykur, and B.D. Park,  {\em Constructing infinitely many smooth structures on small 4-manifolds}, J. Topol. {\bf 1} (2008), 409--428.

\bibitem[2]{AP1} A. Akhmedov and B.D. Park, {\em Exotic smooth structures on small $4$-manifolds}, Invent. Math. {\bf 173} (2008), 209--223.
 
\bibitem[3]{AP2} A. Akhmedov and B.D. Park, {\em Exotic smooth structures on small $4$-manifolds with odd signatures}, Invent. Math. {\bf 181} (2010), 577--603.
 
\bibitem[4]{ADK} D. Auroux, S. Donaldson, and L. Katzarkov, {\em Luttinger surgery along Lagrangian tori and non-isotopy for singular symplectic plane curves}, Math. Ann. {\bf 326} (2003), 185--203.

\bibitem[5]{BK}  S. Baldridge and P. Kirk, {\em A symplectic manifold homeomorphic but not diffeomorphic to $\CPC$}, Geom. Topol. {\bf 12} (2008), 919--940.

\bibitem[6]{RevEng} R. Fintushel, B.D. Park, and R. Stern, {\em Reverse engineering small $4$-manifolds}, Algeb. Geom. Topol. {\bf 7} (2007), 2103Ð-2116. 

\bibitem[7]{rat} R. Fintushel and R. Stern, {\em Rational blowdowns of smooth $4$-manifolds}, Jour. Diff. Geom. {\bf 46} (1997), 181--235.

\bibitem[8]{KL4M} R. Fintushel and R. Stern, {\em Knots, links, and $4$-manifolds}, Invent. Math. {\bf 134} (1998), 363--400.

\bibitem[9]{DN} R. Fintushel and R. Stern, {\em Double node neighborhoods and families of simply connected 4-manifolds with $b^+=1$},  J. Amer. Math. Soc. 19 (2006), 171--180.

\bibitem[10]{Santeria} R. Fintushel and R. Stern, {\em Surgery on nullhomologous tori and simply connected 4-manifolds with $b^+=1$}, J. Topol. {\bf 1} (2008), 1--15.

\bibitem[11]{MS} D. McDuff and M. Symington, {\em Associativity properties of the symplectic sum},
Math. Res. Lett. {\bf 3} (1996), 591--608.

\bibitem[12]{MMS}  J. Morgan, T. Mrowka, and Z. Szab\'{o}, {\em Product formulas along $T^3$ for Seiberg-Witten invariants}, Math. Res. Letters {\bf 4} (1997), 915--929.

\bibitem[13]{OR} P. Orlik and F. Raymond, {\em Actions of the torus on $4$-manifolds. 1}, Trans. Amer. Math. Soc. {\bf 152} (1970), 531--559. 

\bibitem[14]{S} M. Symington, {\em A new symplectic surgery: the 3-fold sum},  Topology Appl. {\bf 88} (1998), 27--53.

\bibitem[15]{T} C. Taubes, {\em More constraints on symplectic manifolds from
Seiberg-Witten invariants}, Math. Res. Letters   {\bf 2} (1995), 9--14.

\end{thebibliography}
\end{document}